\documentclass{amsart}

\usepackage{a4wide}

\usepackage[english]{babel}

\usepackage{amsmath}
\usepackage{amssymb}

\usepackage{amscd}
\input xypic
\usepackage{tabularx}

    \newtheorem{lem}{Lemma}[section]
    \newtheorem{prop}[lem]{Proposition}
    \newtheorem{thm}[lem]{Theorem}

   \theoremstyle{definition}
    \newtheorem{Def}[lem]{Definition}

\newtheorem{List}[lem]{List}

\theoremstyle{remark}
    \newtheorem{rem}[lem]{Remark}

\def\duale{\check{\ }}
\def\la{\langle}
\def\ra{\rangle}

\def\phi{\varphi}

\newcommand{\GL}{\mathrm{GL}}
\newcommand{\SL}{\mathrm{SL}}

\newcommand{\PGL}{\mathrm{PGL}}

\def\Z{\mathbf{Z}}
\def\Q{\mathbf{Q}}
\def\R{\mathbf{R}}
\def\C{\mathbf{C}}

\def\pu{\bullet}

\def\mil{\mspace{-9mu}}

\DeclareMathOperator{\im}{Im}

\newcommand{\Pp}[1]{\mathbf{P}^{#1}}
\def\ba{\big|}

\def\A{\mathcal{A}}
\def\M{\mathcal{M}}
\def\X{\mathcal X}

\newcommand{\B}[2]{B({#2},{#1})}
\newcommand{\F}[2]{F({#2},{#1})}

\newcommand{\tB}[2]{\tilde B({#2},{#1})}

\DeclareMathOperator{\Fil}{Fil}
\DeclareMathOperator{\Conf}{Conf}

\def\op{\mathring}

\newcommand{\coh}[3][\Q]{H^{#2}({#3};{#1})}

\title{Rational cohomology of the moduli space of genus 4 curves}
\author{Orsola~Tommasi}
\email{tommasi@math.kun.nl}
\address{Department of Mathematics, University of Nij\-me\-gen,
Toernooiveld, 6525 ED Nij\-me\-gen, The Netherlands}
\date{29$^\text{th}$ September, 2003.}



\begin{document}

\begin{abstract}
We prove that the Poincar\'e polynomial of the moduli space of smooth genus 4 curves is $1+t^2+t^4+t^5$. 
We show this by producing a stratification of the space, such that all strata are geometric quotients of complements of discriminants.

\end{abstract}

\maketitle
\section{Introduction and results}

In this paper the rational cohomology of the moduli space $\mathcal M_4$ of non-singular complex genus 4 curves is computed. This is achieved by considering a natural stratification of $\mathcal M_4$, and determining the cohomology of each stratum. Non-singular genus 4 curves can be divided into 3 classes (see \cite{Fa}):\begin{enumerate}
\item\label{I} curves whose canonical model is the complete intersection of a cubic surface and a non-singular quadric surface in $\Pp3$;
\item\label{II} curves whose canonical model is the complete intersection of a cubic surface and a quadric cone in $\Pp3$;
\item\label{III} hyperelliptic curves.
\end{enumerate}

Denote by $C_0$ the locus in $M_4$ of curves of type \ref{I}, by $C_1$ the locus of curves of type \ref{II} and by $C_2$ the hyperelliptic locus. We have a three-steps filtration
\begin{equation}\label{filtraz} C_2\subset \overline{C_1}\subset\overline{C_0}=\mathcal M_4.\end{equation}

The space $C_2$ can be studied from the theory of binary forms, and it is easy to show that it has the rational cohomology of a point.
The spaces $C_0$ and $C_1$ are moduli spaces of smooth complete intersections, and their rational cohomology was not known. 
We choose to consider the elements of $C_0$ as representing isomorphism classes of non-singular curves of type $(3,3)$ on $\Pp1\times\Pp1$. Analogously, elements of $C_1$ can be regarded as isomorphism classes of non-singular curves of degree 6 on the weighted projective space $\Pp{}(1,1,2)$. The space $C_0$ is then the quotient of an open subset of the space of bihomogeneous polynomials of bidegree $(3,3)$ in the two sets of indeterminates $x_0,x_1$ and $y_0,y_1$, by the action of the automorphism group of $\Pp1\times\Pp1$. Analogously, $C_1$ is the quotient of an open subset of $\Gamma(\mathcal O_{\Pp{}(1,1,2)}(6))$ by the action of the automorphism group of $\Pp{}(1,1,2)$.

This means that in both cases we are in the situation of
the paper \cite{PS}. As a standard choice of notation, we denote by $V$ the vector space (of dimension $N$), by $X$ the open subset which we are interested in quotienting, and 
by $\Sigma$ the complement of $X$ in $V$, which is a hypersurface called the \emph{discriminant hypersurface}. We denote by $G$ the group acting, by
$\phi: X\rightarrow X/G$ the geometric quotient, and by $\rho$ the orbit inclusion
$$\begin{matrix} G&\longrightarrow&X\\g&\longmapsto&g(x)\end{matrix}$$
where $x$ is a fixed point of $X$.
Recall that $H^\pu(G)$ is an exterior algebra freely generated by classes $\eta_i\in H^{2r_i-1}(G)$. 
We show that in both cases in which we are interested, the following generalization of Leray-Hirsch theorem applies:
\begin{thm}[\cite{PS}]\label{ps}
Suppose there are subschemes $Y_i\subset\Sigma$ of pure codimension $r_i$ in $V$ whose fundamental classes map to a non-zero multiple of $\eta_i$ under the composition 
$$\bar H_{2(N-r_i)}(Y_i)\rightarrow \bar H_{2(N-r_i)}(\Sigma) \xrightarrow{\sim} H^{2r_i-1}(X)\xrightarrow{\rho^*}H^{2r_i-1}(G).$$

Denote the image of $[Y_i]$ in $H^\pu(X;\Q)$ by $y_i$; then the map $a\otimes\eta_i\mapsto\phi^*a \cup y_i$, $a\in H^\pu(X/G;\Q)$ extends to an isomorphism of graded $\Q$-vector spaces 
$$H^\pu(X/G;\Q)\otimes H^\pu(G;\Q)\xrightarrow{\sim}H^\pu(X ;\Q).$$
\end{thm}

We use Theorem \ref{ps} to interpret $H^\pu(X;\Q)$ as a tensor product of the rational cohomology of $X/G$ and $G$. This allows us to recover $H^\pu(C_j;\Q)$ from the other data. 
The problem is then essentially reduced to the computation of the rational cohomology of $X$. Vassiliev invented a general topological method for calculating the cohomology of complements of discriminants, or, more generally, of spaces of non-singular functions. First, the computation of the cohomology of $X$ is reduced to that of the (Alexander dual) Borel-Moore homology of its complement $\Sigma$. Secondly, Vassiliev constructed a simplicial resolution of $\Sigma$, and a stratification of this resolution, based on the study of possible configurations of the singularities of the elements of $\Sigma$. The spectral sequence associated to this filtration is proved to converge to the Borel-Moore homology of $\Sigma$.
Vassiliev used this method in \cite{Vart} to determine, for instance, the real cohomology of the space of non-singular quadric curves in $\Pp2$, and of non-singular cubic surfaces in $\Pp3$. Gorinov (\cite{G}) modified Vassiliev's method so that it applies to a wider range of situations. In this way he could calculate the real cohomology of the space of quintic curves in $\Pp2$. We compute the cohomology of the complements of the discriminants we are interested in, by a modification of Gorinov-Vassiliev's method. We include in this paper an exposition of the method we use. The functoriality of the whole construction is emphasized by describing it in the language of cubical schemes.
Note that Gorinov-Vassiliev's method respects the mixed Hodge structure of the space $X$. Hence we can use it to find the mixed Hodge structure of the cohomology of the complement of the discriminant.

By applying the techniques above, we will show the following.
\begin{thm}\label{c0}
The rational cohomology of the space $C_0$ of non-singular genus 4 curves whose canonical model lies on a non-singular quadric in $\Pp3$ is
$$\coh k{C_0}=\left\{\begin{array}{ll} \Q &k=0\\
\Q(-3)&k=5\\
\,0 &\text{else.}\end{array}\right.$$
\end{thm}

\begin{thm}\label{c1}
The space $C_1$ of smooth genus 4 curves whose canonical model lies on a quadric cone in $\Pp3$ has the rational cohomology of a point.
\end{thm}

Theorems \ref{c0} and \ref{c1} allow us to compute the spectral sequence associated to the filtration (\ref{filtraz}). This establishes the main result:

\begin{thm}\label{main}
The rational cohomology of $\M_4$ is as follows:
$$\coh k{\M_4}=\left\{\begin{array}{ll} \Q & k=0\\\Q(-1) & k=2\\\Q(-2) & k=4 \\\Q(-3) & k=5 \\ \,0 & else.\end{array}\right.$$
\end{thm}

This result agrees with what was previously known about the cohomology of $\mathcal M_4$. In particular, its Euler characteristic had been computed by Harer and Zagier in \cite{HZ}, and found to be 2. 

The plan of the paper is as follows. In Section \ref{ttools} we explain our version of Gorinov-Vassiliev's method for computing cohomology of complements of discriminants. We formulate it by using the language of $\A$-cubical schemes, which is also introduced. Moreover, we show or recall some useful homological results. In Section \ref{C_0} and Section \ref{C_1} we calculate explicitly the rational cohomology of respectively $C_0$ and $C_1$, establishing Theorems \ref{c0} and \ref{c1}. In Section \ref{C_2} we give a short proof, in the style of the paper, of the fact that the moduli space of hyperelliptic curves of genus $g\geq2$ has the rational cohomology of a point.

\subsection{Notations and conventions}
The symbol $\C^n$ will denote the Euclidean space $\C^n\cong\mathbf R^{2n}$.
The complex projective space of dimension $n$ over $\C$ will be denoted by $\Pp n$. The Grassmannian of linear subspaces of dimension $m$ in the vector space $V$ will be denoted by $G(m,V)$.
The symbol $S_n$ denotes the symmetric group in $n$ elements, $\GL(n)$ and $\PGL(n)$ indicate respectively the general linear group and the projective linear group of dimension $n$ over $\C$.

The symbol $\bar H_\pu(Z,S)$ will denote Borel-Moore homology of the space $Z$ with the local system of coefficients $S$. In this paper we make an extensive use of Borel-Moore homology, i.e., homology with locally finite support. A reference for its definition and for the properties we use can be for instance \cite[Chapter 19]{Fulton}.

$\Delta_N$ will denote the $N$-dimensional standard simplex in $\R^{N+1}$.

All complex varieties will be considered with the analytic topology.
We denote by $\Q(k)$ the Tate Hodge structure on $\Q$ of weight $-2k$ (for the definition, see \cite{HII}).

If not otherwise specified, all cohomology and Borel-Moore homology groups are considered with rational coefficients.

\section{Topological tools}\label{ttools}

\subsection{Gorinov-Vassiliev's method}\label{method}

We explain here the method we use for computing the cohomology of complements of discriminants. It has been formulated by Vassiliev (see for instance \cite{Vart}), and successively modified by Gorinov (\cite{G}) in a way that allows it to be applied to a wider range of situations.
We slightly modified Gorinov-Vassiliev's method though, in order to adapt it to configurations on a projective algebraic variety, not necessarily smooth. 

\begin{Def}
Let $Z$ be a projective variety. A subset $S\subset Z$ is called a \emph{configuration} in $Z$ if it is compact and non-empty.
The space of all configurations in $Z$ is denoted by $\Conf(Z)$.
\end{Def}

\begin{prop}[\cite{G}]
The space of all configurations in $Z$ has the structure of a compact complete metric space.
\end{prop}

Let us fix $Z$, and consider a vector space $V$ such that there is a map
$$\begin{matrix}V & \longrightarrow &\Conf(Z)\cup \{\emptyset\}\\
v& \longmapsto & K_v,
\end{matrix}$$

such that $K_0=Z$, and
$$L(K):=\{v\in V: K\subset K_v\}$$ 
is a linear space for all $K\in \Conf(Z)$. 

The most natural case is that in which the elements of $V$ can be seen as sections in a line bundle on $Z$ and $K_v$ is the set of singular points of $v$. 

We define the discriminant as 
$$\Sigma:=\{v\in V: K_v\neq\emptyset\}.$$

The method is based on the fact that there is a direct relation between the cohomology of the complement of $\Sigma$ in $V$ and the Borel-Moore homology of $\Sigma$.
The cap product with the fundamental class $[\Sigma]$ of the discriminant induces for all indices $i$ an isomorphism
$$\tilde H^i(V\setminus\Sigma;\Q)\cong H^{i+1}(V,V\setminus\Sigma;\Q)\overset{\cap[\Sigma]}{\cong}\bar H_{2M-1-i}(\Sigma;\Q)(-M),$$
where we have denoted by $M$ the complex dimension of $V$.

Our aim is to compute the Borel-Moore homology of $\Sigma$. Gorinov constructed a simplicial resolution for $\Sigma$, starting from a collection $X_1,\dots,X_N$ of families of configurations in $Z$. For his construction to work, the $X_i$'s have to satisfy some axioms. We list below the properties which we will require for $X_1,\dots,X_N$.  

\begin{List}\label{ax}\ 
\begin{enumerate}
\item\label{primo} For every element $v\in\Sigma$, $K_v$  must belong to some $X_i$.
\item\label{order} If $x\in X_i$, $y\in X_j$, $x\subsetneq y$, then $i<j$.
\item\label{grass} For every index $j=1,\dots,N$ the space
$L(x)\subset\Sigma$
 has the same dimension $d_j$ for every configuration $x\in X_j$.
Moreover, for all indices $j$ the map:
$$\begin{matrix}\phi_j:& X_j & \longrightarrow & G(d_j,V)\\
& v &\longmapsto & L(K_v)\end{matrix}$$
from the configuration space $X_j$ to the Grassmannian of linear subspaces of dimension $d_j$ of $V$, is continuous.
\item\label{disj} $X_i\cap X_j=\emptyset$ for $i\neq j$.
\item\label{border} Any $x\in \bar X_i\setminus X_i$ belongs to some $X_j$ with $j<i$.
\item\label{loctri} For every $i$ the space 
$$\mathcal T_i=\{(p,x)\in Z\times X_i: p\in x\},$$
with the evident projection,
is the total space of a locally trivial bundle over $X_i$.
\item\label{subconf} Suppose $X_i$ consists of finite configurations. Then for all $y,x$ such that $x\in X_i$, $y\subsetneq x$, the configuration $y$ belongs to $X_j$ for some index $j<i$.  
\end{enumerate}
\end{List}

Note that the maps $\phi_j$ of condition \ref{grass} are always continuous, if 
$L(\{z\})$ has the same dimension $d$ for all $z\in Z$ and the map to the Grassmannian of linear subspaces of dimension $d$ of $V$,

$$\begin{matrix} L:&Z&\longrightarrow& G(d,V)\\
&z&\longmapsto&L(\{z\})\end{matrix}$$
is continuous.

In the construction of the resolution, we will use the language of \emph{cubical spaces}.

\begin{Def}
A \emph{cubical space} over an index set $\A$ (briefly, an $\A$-cubical space) is a collection of topological spaces $\{Y(I)\}_{I\subset\A}$ such that for each inclusion $I\subset J$ we have a natural continuous map $f_{IJ}: Y(J)\rightarrow Y(I)$ such that $f_{IK}=f_{IJ}\circ f_{JK}$ whenever $I\subset J\subset K$. 
\end{Def}

We are ready to define the cubical space and the index set we work with.

$$\mathcal A=\{1,2,\dots,N\}$$

$$\Lambda(I):=\{A\in\prod_{i\in I} X_i:i<j\Rightarrow A_i\subset A_j\}.$$

$$\X(I):=\{(F,A)\in\Sigma\times\Lambda(I):K_F\supset A_{\max(I)}\}\text{ if }I\neq\emptyset,$$

$$\X(\emptyset):=\Sigma.$$

Analogously, we define the following (auxiliary) cubical spaces:

$$\tilde\Lambda(I):=\{A\in\prod_{i\in I}\overline{X}_i: i<j\Rightarrow A_i\subset A_j\},$$
$$\tilde\X(I):=\{(F,A)\in \Sigma\times\tilde\Lambda(I):K_F\supset A_{\max(I)}\}\text{ if }I\neq\emptyset,$$
$$\tilde\X(\emptyset):=\Sigma.$$

The subsets $I\subset\A$ represent the vertices of the cube $\square_\A$, which are in one to one correspondence with the faces $\Delta_I$ of the simplex
$$\Delta_\A:=\left\{(f:\A\rightarrow[0,1]):\sum_{a\in\A}f(a)=1\right\}.$$
By definition, $\Delta_I=\{f\in\Delta_\A:f|_{\A-I}=0\}$. 
To every inclusion $I\subset J$ there is a natural map associated, which we denote by $e_{IJ}:\Delta_I\rightarrow\Delta_J$.

\begin{Def}
Let $Y(\pu)$ be a cubical space over an index set $\A$.

Note that $Y(\pu)$ has a natural augmentation towards $Y(\emptyset)$. 
Then the \emph{geometric realization} of $Y(\pu)$ is defined as the map
$$|\epsilon|: \ba Y(\pu)\ba \longrightarrow Y(\emptyset)$$
induced from the natural augmentation on the space
$$\ba Y(\pu)\ba =\coprod_{I\subset\A}(\Delta_I\times Y(I))/R,$$
where $R$ is the equivalence relation given by 
$$(f,y)\ R\ (f',y')\Leftrightarrow f'=e_{IJ}(f),y=f_{IJ}(y').$$
\end{Def}

Note that $\Delta_\emptyset=\emptyset$ so that $Y(\emptyset)$ does not appear in the construction of the geometric realization. 

We construct the geometric resolution of all the cubical spaces defined above.
We can define a surjective map $\phi$ between $\ba\tilde\Lambda(\pu)\ba$
 to $\ba\Lambda(\pu)\ba$ as follows. Let $(t,A)\in\Delta_I\times\tilde\Lambda(I)$, and $[t,A]$ the corresponding class in $\ba\tilde\Lambda(\pu)\ba$. 
Note that by Conditions 4 and 5 in List~\ref{ax}, for each $A_i$ ($i\in I$) there exists a unique index $k(i)\in\A$ such that $A_i\in X_{k(i)}$.

We define $\phi([t,A])$ as the class in $\ba\Lambda(\pu)\ba$ of the element $(s,B)\in \Delta_J\times\Lambda(J)$, where 

$$J:=\{k\in\A: A_i\in X_k \text{ for some }i\in I\}=\{k(i):i\in I\};$$
$$B:=\prod_{k\in J}B_k,\ B_k:=A_i\text{ for any index }i: k(i)=k;$$
$$s:J\longrightarrow [0,1],\ s(k):=\sum_{i\in I:k(i)=k}t(i).$$

Observe that $\phi$ acts by contracting all closed simplices in $\ba\tilde\Lambda(\pu)\ba$ corresponding to inclusions of configurations of the form $x=x=x=\dots=x$.

We define analogously the map
$$\begin{matrix}\psi: &\ba\tilde\X(\pu)\ba&\longrightarrow &\ba\X(\pu)\ba\\
&[(F,A),t]&\longmapsto&(F,\phi(A,t)).\end{matrix}$$

In the rest of the paper, we consider the spaces $\ba\tilde\Lambda(\pu)\ba$ and $\ba\tilde\X(\pu)\ba$ with the quotient topology under the equivalence relation $R$ of the direct product topology of the $\tilde\Lambda(I)$'s (respectively, $\X(I)$'s).
The topology on $\ba\Lambda(\pu)\ba$ and $\ba \X(\pu)\ba$ is the topology induced by $\phi$ (respectively, $\psi$).


\begin{prop}[\cite{G}]
The geometric realization of $\X(\pu)$,
$$|\epsilon|: \ba \X(\pu)\ba \longrightarrow \X(\emptyset)=\Sigma,$$
is a homotopy equivalence and induces an isomorphism on Borel-Moore homology groups.
\end{prop}

\proof
It is enough to prove that $|\epsilon|$ is a proper map, and that its fibers are contractible. 

Fix a compact subset $W\subset\Sigma$. We define a cubical space as follows:
$$\tilde\X_W(I):=\{(F,A)\in W\times\tilde\Lambda(I):K_F\supset A_{\max(I)}\}\text{ if }I\neq\emptyset,$$
$$\tilde\X_W(\emptyset):=W.$$

Note that $\tilde\X_W(I)$ is compact for all $I\in\A$, and so is the space $\ba\tilde\X_W(\pu)\ba$.

Then $|\epsilon|^{-1}(W)$ is compact, because it coincide with the image of the continuous map
$$\begin{matrix}\psi_W: &\ba\tilde\X_W(\pu)\ba&\longrightarrow &\ba\X(\pu)\ba\\
&[(F,A),t]&\longrightarrow&(F,\phi(A,t)).\end{matrix}$$
Hence the map $\ba\epsilon\ba$ is proper.

We show next that the fibres of $\ba\epsilon\ba$ are contractible.
Consider the fiber over $v\in\Sigma$. 
By conditions \ref{primo} and \ref{disj} of List \ref{ax}, there is a unique index $j$ such that $K_v\in X_{j}$.
 By definition, 
$|\epsilon|^{-1}(v)$ is a cone with vertex $[(K_v,v)]\in \ba \X|_{\{j\}}(\pu)\ba \hookrightarrow \ba \X(\pu)\ba $, so it is clearly contractible.  
\qed

\medskip
For every index set $I\subset\A$, we can restrict $\Lambda(\pu)$ and $\X(\pu)$ to the index set $I$, getting the two $I$-cubical spaces $\Lambda|_I(\pu)$ and $\X|_I(\pu)$. Then for every $I\subset J\subset\A$, there are natural embeddings $\ba \Lambda|_I(\pu)\ba \hookrightarrow\ba \Lambda|_J(\pu)\ba $ and $\ba \X|_I(\pu)\ba \hookrightarrow\ba \X|_J(\pu)\ba $. 
In this way we can define an increasing filtration on $\ba\Lambda(\pu)\ba $ by posing
$$\Fil_j\ba \Lambda(\pu)\ba :=\im\left(\ba \Lambda|_{\{1,2,\dots,j\}}\ba \hookrightarrow\ba \Lambda(\pu)\ba \right)$$ for $j=1,\dots,N$. We define analogously the filtration $\Fil_j\ba \X(\pu)\ba $ on $\ba \X(\pu)\ba $. We use the notation $F_j:=\Fil_j\ba \X(\pu)\ba\setminus\Fil_{j-1}\ba \X(\pu)\ba $, $\Phi_j:=\Fil_j\ba \Lambda(\pu)\ba\setminus\Fil_{j-1}\ba \Lambda(\pu)\ba $.

The filtration $\Fil_j\ba \X(\pu)\ba $ defines a spectral sequence that converges to the Borel-Moore homology of $\Sigma$. Its term $E^1_{p,q}$ is isomorphic to $\bar H_{p+q}(F_p;\Q)$.

\begin{prop}[\cite{G}]\label{ucci}
\renewcommand{\labelenumi}{\arabic{enumi}.}
\begin{enumerate}
\item For every $j=1,\dots,N$, the stratum $F_j$ is a complex vector bundle of rank $d_j$ over $\Phi_j$.
The space $\Phi_j$ is in turn a fiber bundle over the configuration space $X_j$. 
\item If $X_j$ consists of configurations of $m$ points, the fiber of $\Phi_j$ over any $x\in X_j$ is an $(m-1)$-dimensional open simplex, which changes its orientation under the homotopy class of a loop in $X_j$ interchanging a pair of points in $x_j$.
\item\label{opeco}  If $X_N=\{Z\}$, $F_N$ is the open cone with vertex a point (corresponding to the configuration $Z$), 
over $\Fil_{N-1}\ba \Lambda(\pu)\ba $. 
\end{enumerate}\end{prop}

We recall here the topological definition of an open cone.
\begin{Def}
Let $B$ be a topological space. Then a space is said to be an \emph{open cone} over $B$ with vertex a point if it is homeomorphic to the space 
$$B\times[0,1)/R,$$
where the equivalence relation is $R=(B\times\{0\})^2$.
\end{Def}

\begin{proof}[Proof of Proposition \ref{ucci}]
The second and the third point are clear by construction.
The first point is trivial for $F_j$. For the map $\Phi_j\longrightarrow X_j$, the fiber over a configuration $x$ is given by a union of simplices with vertices determined by the points of $x\subset Z$. Thus the fibration is locally trivial as a consequence of condition \ref{loctri} in List \ref{ax}. 
\end{proof}

\subsection{Homological lemmas}\label{lemmas}

The fiber bundle $(\Phi_j\setminus\Phi_{j-1})\rightarrow X_j$ of Proposition \ref{ucci} is in general non-orientable. As a consequence, we have to consider the homology of $X_j$ with coefficients not in $\Q$, but in some local system of rank one. Therefore we recall here some results and constructions about Borel-Moore homology of configuration spaces, with twisted coefficients.

\begin{Def}
Let $Z$ be a topological space. Then for every $k\geq 1$ we have the space of ordered configurations of $k$ points in $Z$,
$$\F kZ=Z^k\setminus\bigcup_{1\leq i<j\leq k}\{(z_1,\dots,z_k)\in Z^k: z_i=z_j\}.$$
There is a natural action of the symmetric group $S_k$ on $F(k,Z)$. The quotient is called the space of unordered configurations of $k$ points in $Z$,
$$\B kZ=\F kZ/S_k.$$
\end{Def}

The \emph{sign representation} $\pi_1(\B kZ)\rightarrow Aut(\Z)$ maps the paths in $\B kZ$ defining odd (respectively, even) permutations of $k$ points to multiplication by $-1$ (respectively, 1). The local system $\pm\Q$ over $\B kZ$ is the one locally isomorphic to $\Q$, but with monodromy representation equal to the sign representation of $\pi_1(\B kZ)$. We will often call $\bar H_\pu (\B kZ,\pm\Q)$ the \emph{Borel-Moore homology of $\B kZ$ with twisted coefficients}, or, simply, the \emph{twisted Borel-Moore homology of $\B kZ$}. 

\begin{lem}[\cite{Vbook}]\label{lem2}
$\bar H_\pu(\B k{\C^N};\pm \Q)=0$ if $k\geq2$.
\end{lem}

\proof
 In \cite[Theorem 4.3, Corollary 2]{Vbook}, it is proved that the Borel-Moore homology of $\B k{\R^{2N}}$ with coefficients in the system $\pm\Z$ is a finite group. Since $\pm\Q=\pm\Z\otimes_\Z\Q$, the claim follows. 
\qed

\begin{lem}[\cite{Vart}]\label{lem1}
$\bar H_\pu (\B k{\Pp N};\pm\Q)= H_{\pu-k(k-1)}( G(k,\C^{N+1});\Q)$.\\
In particular, $\bar H_\pu (\B k{\Pp N};\pm\Q)$ is trivial if $k>N+1$.
\end{lem}

\begin{lem}
The Poincar\'e polynomial of $\bar H_\pu(\B 2{\C^*});\pm\Q)$ is $t^2(1+t)$. For constant rational coefficients, it is $t^2(1+t)^2$.
\end{lem}

\proof
We will recover the Borel-Moore homology of $\B 2{\C^*}$ from the known situation for $\B 2\C$. Since $\bar H_\pu(\B 2\C;\pm\Q)$ is trivial, we have 
$$\bar H_j(\B2\C;\Q)\cong\bar H_j(\C^2\setminus\{(x,y)\in\C^2:x=y\};\Q)=\left\{\begin{array}{ll}\Q(1) &j=3\\\Q(2) &j=4\\\,0&\text{else}.\end{array}\right.$$

Let us consider the two-steps filtration $\{\{0,a\}:a\in\C^*\}\subset\B2\C$. Clearly it gives a decomposition $\B2\C)\cong\C^*\amalg\B2{\C^*}$. Hence, in the spectral sequence calculating the Borel-Moore homology (respectively the twisted Borel-Moore homology) of $\B2\C$ associated to this filtration we have 
$$E^1_{1,q}=\bar H_{q-1}(\C^*;\Q), $$ $$E^1_{2,q}=\bar H_{q-2}(\B2{\C^*};\Q) \text{ (respectively, } E^1_{2,q}=\bar H_{q-2}(\B2{\C^*};\pm\Q).$$ The Borel-Moore homology of $\C^*$ and of $\B2\C$ being known, this determines the Borel-Moore homology of $\B 2{\C^*}$ both for constant and for twisted rational coefficients.
\qed

\begin{lem}\label{lemQ}
The Poincar\'e polynomial of 
$\bar H_\pu(\B 1{{\Pp1}\times{\Pp1}};\pm\Q)$ is 
$(t^2+1)^2$.\\
The Poincar\'e polynomial of 
$\bar H_\pu(\B 2{{\Pp1}\times{\Pp1}};\pm\Q)$
is $2t^2(t^4+t^2+1)$.\\
The Poincar\'e polynomial of 
$\bar H_\pu(\B 3{{\Pp1}\times{\Pp1}};\pm\Q)$
is $t^4(t^2+1)^2$.\\
The Poincar\'e polynomial of 
$\bar H_\pu(\B 4{{\Pp1}\times{\Pp1}};\pm\Q)$ is $t^8$.\\
The twisted Borel-Moore homology of $\B k{{\Pp1}\times{\Pp1}})$ is trivial for $k\geq5$.
For each of the groups above, the Borel-Moore homology in degree $d$ has Hodge weight $-d$.
\end{lem}

\proof 
We modify here Vassiliev's arguments in the proof of Lemma \ref{lem1} in \cite{Vart}.
The technique we use is that of decomposing $\B k{{\Pp1}\times{\Pp1}}$ into spaces of which the twisted Borel-Moore homology is known. In particular, it is possible to decompose $\Pp 1\times{\Pp1}$ by fixing two lines $l,m$ in different rulings and considering the filtration $S_1\subset S_2\subset S_3\subset S_4$, where
$$S_1:=l\cap m, S_2:=l , S_3:=l\cup m , S_4:={\Pp1}\times{\Pp1}.$$
This means ${\Pp1}\times{\Pp1}$ is the disjoint union of spaces isomorphic to $\{*\},\C,\C,\C^2$ respectively. 

Let us fix $k\geq1$. To any configuration of points in $\B k{{\Pp1}\times{\Pp1}})$ we can associate an ordered partition $(a_1,a_2,a_3,a_4)$, where $a_i$ is the number of points contained in $S_i\setminus S_{i-1}$. We can consider each possible partition of $k$ as defining a stratum in $\B k{{\Pp1}\times{\Pp1}})$, and order such strata by lexicographic order of the index of the partition. Note that all strata with $a_i\geq2$ for some $i$ have no twisted Borel-Moore homology by Lemma \ref{lem2}, so we need not consider them. This is the case for all strata, when $k\geq5$.

As an example, we consider explicitly the case $k=2$. The situation is analogous for the other values. The admissible partitions for $k=2$ are: \begin{description}
\item[(1,1,0,0)] the stratum is isomorphic to $\C$, hence it has cohomology $\Q(1)$ in degree $2$, and trivial elsewhere.
\item[(1,0,1,0)] the stratum is isomorphic to $\C$, hence it has cohomology $\Q(1)$ in degree $2$, and trivial elsewhere.
\item[(1,0,0,1)] the stratum is isomorphic to $\C^2$, hence it has cohomology $\Q(2)$ in degree $4$, and trivial elsewhere.
\item[(0,1,1,0)] the stratum is isomorphic to $\C^2$, hence it has cohomology $\Q(2)$ in degree $4$, and trivial elsewhere.
\item[(0,1,0,1)] the stratum is isomorphic to $\C^3$, hence it has cohomology $\Q(3)$ in degree $6$, and trivial elsewhere.
\item[(0,0,1,1)] the stratum is isomorphic to $\C^3$, hence it has cohomology $\Q(3)$ in degree $6$, and trivial elsewhere.
\end{description}

This gives precisely that $\bar H_j (\B 2{{\Pp1}\times{\Pp1}};\pm\Q)$ is $\Q(j/2)^2$ for $j=2,4,6$, and is trivial otherwise.
\qed

\begin{lem}\label{conticono}
Let us consider a quadric cone $U\subset\Pp3$. Define $U_0$ as $U$ minus its vertex. Then 
$\bar H_\pu(\B 2{U_0};\pm\Q)$ 
has dimension 1 in degree 6 and is trivial otherwise.
The twisted Borel-Moore homology $\bar H_\pu(\B k{U_0};\pm\Q)$ is always trivial for $k\geq3$.
The Hodge weight of the homology in degree $d$ is $-d$.
\end{lem}

\proof Analogous as that of the previous lemma, by using the fact that $U_0$ can be decomposed as the disjoint union of $\C$ (a line of the ruling) and a space homeomorphic to $\C^2$.\qed

\begin{lem}\label{lem:T}
Let us consider the the local system $T$ on $\C^*$ locally isomorphic to $\Q$ and changing its sign when the point moves along a loop in $\C^*$ whose homotopy class is an odd multiple of the generator of $\pi_1(\C^*)\cong\Z$. Then
$\bar H_\pu(\C^*;T)=0$.
\end{lem}

\proof
Let us consider the map $$\begin{array}{lccc}q: &\C^*&\longrightarrow&\C^*\\
&t&\longmapsto&t^2.\end{array}$$
Clearly $q_*\Q = \Q\oplus T$. Since $\bar H_\pu(\C^*;q_*\Q)\cong \bar H_\pu(\C^* ;\Q)$, we have $\bar H_\pu(\C^*;T)=0$.
\qed

\subsection{Some cases where the Borel-Moore homology is trivial}\label{simplices}

In Sections \ref{C_0} and \ref{C_1} we use the fact that most strata of the filtration we consider over the geometric realization of discriminants give no contribution to its Borel-Moore homology.
We give here some ideas about the reason why it happens.
In general they are either situations where the space is a locally trivial fiber bundle whose fibers have trivial Borel-Moore homology in the induced system of coefficients, or situations were non-discrete configurations are involved.

\begin{lem}\label{opencone}
Let $C$ be an open cone with vertex a point over a space $B$. Then there are the following isomorphisms:
$$\bar H_\pu(C;\Q)\cong H_\pu(C,B;\Q)\cong H_{\pu-1}(B,\text{point};\Q).$$
\end{lem}

\proof
The first isomorphism comes from the characterization of Borel-Moore homology as relative homology of the one point compactification of the space modulo the added point. The second is the border isomorphism of the exact sequence associated to the triple $(\text{point},B,C)$. 
\qed

\begin{lem}\label{lemline} Suppose we have a variety $Z$ and the following families of configurations in $Z$:
$$\begin{array}{l} X_1=\B 1Z;\\
X_2=\{\{p,q\}\in \B2Z:\text{ $p$ and $q$ lie on a line }l\subset Z\};\\
X_3=\{\{p,q,r\}\in \B3Z:\text{ $p,q$ and $r$ lie on a line }l\subset Z\};\\
X_4=\{\text{lines on }Z\}.
\end{array}$$
Construct the cubical space $\Lambda(\pu)$, its geometric realization and the filtration as in Section \ref{method}. Then the space $\Phi_4$ has trivial Borel-Moore homology.
\end{lem}

\proof
The space $\Phi_4$ is a fiber bundle over $X_4$. Its fibre $\Psi$ over a line $l\subset Z$ is the union of all simplices with 4 vertices $p,q,r,l$ for all $\{p,q,r\}\in\B3Z$. Note that we have to take \emph{closed} simplices minus their face with vertices $p,q,r$, so they are indeed open cones with vertex $l$ over the closed simplices with vertices corresponding to the points $p,q,r$.
The system of coefficients induced on the fibre changes its sign if we interchange two of the points $p,q,r$. 
The union of all \emph{open} simplices with vertices $p,q,r,l$ is a non-oriented simplices bundle over $\B3l$. Hence the Borel-Moore homology of the union of open simplices with three vertices on $l$ is trivial, because $\bar H_\pu(\B3l;\pm\Q)$ is trivial (see Lemma \ref{lem1}).
This means that we need to consider only the Borel-Moore homology of the union of the external faces of the simplices considered before. This space admits a filtration of the form: $A_0=\{l\}$, $A_1$ is the union of the open segments joining $l$ and a point of $l$, $A_2$ is the union of the open simplices with vertices $l$ and two distinct points on $l$. This gives a spectral sequence converging to the Borel-Moore homology of $\Psi$, with the following $E^1$-term:
$$\begin{array}{r|ccc}
 1&&\Q&\Q\\
 0&&&\\
-1&\Q&\Q&\\
\hline
&1&2&3
\end{array}$$
Our space is an open cone, and it is a consequence of Lemma \ref{opencone} that open cones have trivial Borel-Moore homology in degree 0. This implies that the row $q=-1$ of the spectral sequence is exact.
As for the row $q=1$, it follows from the shape of a generator of $E^1_{2,1}$ that it maps to a generator of $E^1_{1,1}$ under the differential $d^1$. Then the Borel-Moore homology of $\Psi$ is trivial, which proves the claim.
\qed
\begin{rem}\label{remlp}
A slight modification of Lemma \ref{lemline} allows us to conclude that also strata which singular configurations which are union of a line and a fixed finite number $k$ of points have trivial Borel-Moore homology. 

As before, we consider the projection to the configuration space and look at the fiber $\Psi$ over a fixed configuration $\{a_1,a_2,\dots a_k\}\cup l$. 
Let $h\geq3$ be the maximal number of isolated points lying on the same line which appear in the previous strata of $\ba \Lambda(\pu)\ba $. 
For any choice of distinct points $b_1,\dots,b_h$ on $l$, the union of the $(k+h)$-dimensional open simplices with vertices identified with $a_1,\dots,a_k,b_1,\dots,b_h,l$ is contained in $\Psi$. Indeed, $\Psi$ is the union of the open simplex $D$ with vertices $a_1,\dots,a_k,l$ and the union of such simplices for every choice of $b_1,\dots,b_h$. This means that $\Psi$ has a natural projection $\pi$  to $D$. Each point $p$ of $\Psi\setminus D$ is contained in exactly one open simplex with vertices $b_1,\dots,b_h,d$ for some $d\in D$; we pose $\pi(p)=d$. On $D\subset\Psi$ the projection $\pi$ coincides with the identity. If we look at the fibres of $\pi$, we see that they are homeomorphic to the fibres of the map $\Phi_4\longrightarrow X_4$ studied in the proof of Lemma \ref{lemline}. Then we can apply the result found there, and conclude that the fibres of $\pi$ have trivial Borel-Moore homology. Then also $\Psi$ has trivial Borel-Moore homology, which implies the claim.\end{rem}

We consider next the case of the union of two rational curves, intersecting in one point. For simplicity, we state the result in the case of lines.

\begin{lem}\label{2lines} Suppose we have a variety $Z$ and the following families of configurations in $Z$:
$$\begin{array}{l} X_1=\B 1Z;\\
X_2=\B2Z;\\
X_3=\B3Z;\\
X_4=\B4Z;\\
X_5=\{\text{lines on }Z\};\\
X_6=\{\{p,q,r,s,t\}\in \B5Z:\text{ $p,q$ lie on a line }l\subset Z,\\\text{ $r,s$ lie on a line }m\subset Z, l\cap m=\{t\}\};\\
X_7=\{l\cup\{p\}: l\subset Z\text{ line}, p\notin l\};\\
X_8=\{l\cup\{p,q\}: l\subset Z\text{ line}, p,q\notin l, p\neq q\};\\
X_{9}=\{l\cup m: l,m\subset Z\text{ line}, \#(l\cap m)=1\}.
\end{array}$$
Construct the cubical space $\Lambda(\pu)$, its geometric realization and the filtration as in Section \ref{method}. Then the space $\Phi_{9}$ has trivial Borel-Moore homology.
\end{lem}

\proof
Consider the fiber $\Psi$ of the projection $\Phi_8\longrightarrow X_8$ over a configuration $l\cup m$, such that $l\cap m=\{t\}$. The maximal chains of inclusions of configurations we can construct are of the following forms:
\renewcommand{\labelenumi}{(\alph{enumi})}
\begin{enumerate}
\item$\{p_1,p_2,p_3,p_4\}\subset l \subset l\cup\{q_1,q_2\}\subset l\cup m$, where $p_i\in l$, $q_j\in m\setminus \{t\}$;
\item$\{q_1,q_2,q_3,q_4\}\subset m \subset m\cup\{p_1,p_2\}\subset l\cup m$, where $q_i\in m$, $p_j\in l\setminus \{t\}$;
\item$\{p_1,p_2,t,q_1,q_2\}\subset l\cup\{q_1,q_2\}\subset l\cup m$, where $p_i\in l\setminus \{t\}$, $q_j\in m\setminus \{t\}$;
\item$\{p_1,p_2,t,q_1,q_2\}\subset m\cup\{p_1,p_2\}\subset l\cup m$, where $p_i\in l\setminus \{t\}$, $q_j\in m\setminus \{t\}$.
\end{enumerate}
The simplices constructed from chains of inclusions of type (c) and (d) are contained in the simplices arising from chains of type (a) and (b), so that we have to consider just these. 
Denote by $U_1$ the union of closed simplices with vertices given by the configurations $p_1,p_2,p_3,p_4,q_1,q_2,l$ with $p_i\in l$ and $q_j\in m\setminus \{t\}$. Analogously, denote by $U_2$ the union of closed simplices with vertices given by the configurations $q_1,q_2,q_3,q_4,p_1,p_2,m$. 
The fiber $\Psi$ is then the open cone (with vertex corresponding to the configuration $l\cup m$) over $U_1\cup U_2$. Note that the intersection of $U_1$ and $U_2$ is given by the union of closed simplices with vertices $p_1,p_2,q_1,q_2,t$. 
It follows from the results of Remark \ref{remlp} that the Borel-Moore homology of the open cone over $U_1\cap U_2$ is trivial. The same can be said for the open cone over the union of simplices with 4 points of $l$ and two points of $m$ as vertices. The remaining part of $U_1$ is an open cone with vertex $l$ (over simplices with 4 points of $l$ and two points of $m$ as vertices). This space is clearly contractible, so that the open cone above it is contractible by Lemma \ref{opencone}.  We can conclude exactly the same for $U_2$. Then the claim holds. 
\qed

\begin{lem}
Suppose $Z$ is the product of $\C$ and a variety $M$. Consider the following families of configurations in $Z$:
$$\begin{array}{l}
X_1=\B1Z;\\
X_2=\B2Z;\\
X_3=\{\{a,b,c\}\in \B3Z:a,b,c\in \C\times\{p\}, p\in M\};\\
X_4=\{\{a,b,c\}\in \B3Z:b,c\in \C\times\{p\}, p\in M,a\notin \C\times\{p\}\};\\
X_5=\{\{a,b,c,d\}\in \B4Z:a,b,c,d\in \C\times\{p\}\};\\
X_6=\{\{a,b,c,d\}\in \B4Z:a,b\in \C\times\{p\},c,d\in \C\times\{q\},p\neq q\}.
\end{array}$$
Construct the cubical space $\Lambda(\pu)$, its geometric realization and the filtration as in Section \ref{method}. Then the space $\Phi_6$ has trivial Borel-Moore homology.
\end{lem}

\proof
The space $\Phi_6$ is an open non-orientable simplicial bundle over $X_6$. We study the Borel-Moore homology of $X_6$ in the system of coefficients locally isomorphic to $\Q$, with orientation induced by the orientation of the simplices. We look at the ordered situation. Let $Y=\F2M\times\F2\C\times\F2\C$. Every point $(p,q,a,b,c,d)\in Y$ gives an ordered configuration of points $((a,p),(b,p),(c,q),(d,q))$, and the twisted Borel-Moore homology of $X_6$ can be identified with the part of $\bar H_\pu(Y;\Q)$ which is \begin{itemize}
\item[-] anti-invariant under the action of loops interchanging $a$ and $b$;
\item[-] anti-invariant under the action of loops loops interchanging $c$ and $d$;
\item[-] invariant under the action of loops interchanging $p$ and $q$, $a$ and $c$, $b$ and $d$. 
\end{itemize}
It is clear that such homology classes cannot exist, because $$\bar H_\pu(Y;\Q)=\bar H_\pu(\F2M;\Q)\otimes\bar H_\pu(\F2\C;\Q)\otimes \bar H_\pu(\F2\C;\Q)$$ and $\bar H_\pu(\F2\C;\Q)$ contains no classes that are anti-invariant with respect to the interchange of points (see Lemma \ref{lem2}).
\qed

A variation of the situation of the above lemma is given, for instance, by the case of configurations of two triplets of collinear points in $\Pp1\times\Pp1$. In that case we have to use the fact that there are no anti-invariant homological classes in $\bar H_\pu(\F3{\Pp1};\Q)$.

\section{Curves on a non-singular quadric}\label{C_0}

Any non-singular quadric surface is isomorphic to the Segre embedding of $\Pp1\times\Pp1$ in $\Pp3$. 
Such a surface is covered by two families of lines: The family of lines of the form $\Pp1\times\{q\}$ (which we call first ruling of the quadric) and that of lines of the form 
$\{p\}\times\Pp1$ (second ruling). A curve on the Segre quadric is always given by the vanishing of a bihomogeneous polynomial in the two sets of variables $x_0,x_1$ and $y_0,y_1$. The bidegree $(n,m)$ of the polynomial has a geometrical interpretation as giving the number of points (counted with multiplicity) in the intersection of the curve with a general line of respectively the second and the first ruling. A curve on $\Pp1\times\Pp1$ is said to be \emph{of type} $(n,m)$ if it is defined by the vanishing of a polynomial of bidegree $(n,m)$.

The curves which are the intersection of the quadric with a cubic surface are the curves of type $(3,3)$. This suggests that the space $C_0$ can be obtained as a quotient of the space of polynomials of bidegree $(3,3)$ by the action of the automorphism group of $\Pp1\times\Pp1$. We denote that vector space by
$$V:=\C[x_0,x_1,y_0,y_1]_{3,3}\cong\C^{16}.$$

In $V$ we can consider the discriminant locus $\Sigma$ of polynomials defining singular curves in $\Pp 1 \times\Pp 1$. The discriminant $\Sigma$ is closed in the Zariski topology, and is a 15-dimensional cone with the origin as vertex. It is also irreducible, because it is the cone over the dual variety of the Segre embedding of the product of two rational normal curves of degree 3. In other words, the projectivization of $\Sigma$ is the dual of the image of the map
$$\xymatrix@R=15pt{{\Pp1\times\Pp1}\ar[r]^{v_3\times v_3}&{\Pp3\times\Pp3}\ar[r]^{\sigma}&{\Pp{}(V)}\\
{([x_0,x_1],[y_0,y_1])}\ar@{|->}[rr]&&[x_0^3y_0^3,x_0^2x_1y_0^3,\dots,x_1^3y_1^3].}$$

The locus of polynomials giving non-singular curves will be denoted by $X=V\setminus \Sigma$.

The automorphism group of $\Pp1\times\Pp1$ is of the form $G\times S_2$. Its torsion part is generated by the involution $\upsilon$ interchanging the two rulings:
$$\begin{matrix}\upsilon:&\Pp1\times\Pp1&\longrightarrow&\Pp1\times\Pp1\\
&([x_0,x_1],[y_0,y_1])&\longmapsto&([y_0,y_1],[x_0,x_1]).\end{matrix}$$
The connected component containing the identity of the automorphism group is the group
$$G:=\GL(2)\times\GL(2)/\{\lambda I, \lambda^{-1}I\}_{\lambda\in\C^*},$$
and its action on $\Pp1\times\Pp1$ is induced by the action of $\GL(2)$ on $\Pp1$. 

Being the quotient of a reductive group, $G$ is itself reductive. The action of $G$ on binary polynomials induces an action on $V$. The space $X$ is clearly invariant under this action of $G$, and all its point are $G$-stable. Thus there exists a GIT quotient $X/G$, which is a double cover of $C_0$, the involution of the double cover being induced by $\upsilon$.

\subsection{Generalized Leray-Hirsch theorem}\label{LH_0}

The aim of this section is to prove that the action of $G$ on $X$ satisfies the hypotheses of Theorem \ref{ps}. 

We need first to compute the cohomology of $G$.
Consider the map:

$$\begin{array}{cccc}
\iota:& \C^*\times\SL(2)\times\SL(2)&\longrightarrow &G\\
& (\lambda,A,B)&\longmapsto& \left[(\lambda A,\lambda B)\right].
\end{array}$$

The map $\iota$ is an isogeny of connected algebraic groups. Its kernel is finite, hence $\iota$ induces an isomorphism of rational cohomology groups. As a consequence, the cohomology of $G$ is an exterior algebra on three independent generators: a generator $\xi$ of degree 1 and two generators $\eta_1,\eta_2$ of degree 3.

What we need, is to show the surjectivity of the map (induced by the {orbit inclusion} $\rho: G \rightarrow X$)
$$\rho^*: \coh iX\longrightarrow \coh iG.$$

We will do it by studying the map
$$\xymatrix@R=15pt{
*{\psi:}&*{\coh i{X}}\ar@{>}[r]\ar@{<->}[d]&{\coh iG}\ar[r]&{\coh i{\GL(2)\times\GL(2)}}\ar@{<->}[d]
\\&*{\bar H_{31-i}(\Sigma;\Q)}\ar[rr]&&{\bar H_{15-i}((D\times M)\cup(M\times D);\Q),}
}$$
where we have considered the embedding of $\GL(2)$ in the space $M$ of $2\times2$ matrices, and written $D=M\setminus \GL(2)$ for the hypersurface defined by the vanishing of the determinant.
Note that $\coh3G$ is isomorphic to $\coh3{\GL(2)\times\GL(2)}$ and that $\coh1G$ can be identified with the part of $\coh1{\GL(2)\times\GL(2)}$ that is invariant with respect to the interchange of the two factors of the product.

We need only to show that the generators of $\coh\pu G$ are contained in the image of $\psi$. We will see that in this case we can write the generators of the cohomology groups of $X$ quite explicitly, by means of fundamental classes. 

The cohomology of $\GL(2)\times\GL(2)$ is determined by the Borel-Moore homology of the discriminant $D\subset M$. We can compute it by looking at the desingularization
$$\begin{CD} \tilde D=\{(p,A)\in\Pp1\times M: Ap=0\}@>>> D\\
@V\tau VV\\
\Pp1.\end{CD}$$

As a consequence, $\bar H_\pu(\tilde D;\Q) \cong\bar H_{\pu-4}(\Pp1;\Q)$. Hence the only non-trivial groups are $\bar H_6(\tilde D;\Q)$, which is generated by the fundamental class of $\tilde D$, and $\bar H_4(\tilde D;\Q)$, which is generated by the fundamental class of the preimage $\tilde R$ of a point in $\Pp1$. 
As Borel-Moore homology is covariant for proper morphisms, there is a natural map $\bar H_\pu(\tilde D;\Q)\rightarrow\bar H_\pu(D;\Q)$, which must be an isomorphism in degrees 4,6 because in those cases the two groups have the same dimension. Thus we know  generators for $\bar H_\pu(D;\Q)$. The fundamental class of $D$ is a generator of degree 6, and the fundamental class of the image $R$ of $\tilde R$ is a generator of degree 4.
In particular, we can choose $R$ to be the subvariety of matrices with only zeroes on the first column.

We have natural projections
$$D\times M\xrightarrow{a_1} D\xleftarrow{a_2} M\times D$$

and natural immersions

$$\xymatrix@R=0pt{ {D\times M}\ar@{<-}[r]^{i_1} & {D}\ar[r]^{i_2}& {M\times D}\\
{(A,I)}\ar@{<-|}[r]&{A}\ar@{|->}[r]&{(I,A).}}$$

Then we have 
$$\bar H_{14}((D\times M)\cup (M\times D);\Q)\cong \bar H_{14}(D\times M;\Q)\oplus\bar H_{14}(M\times D;\Q) \cong \Q\la a_1^*([D]),a_2^*([D])\ra,$$
and the part which comes from the cohomology of $G$ is generated by $a_1^*([D])+a_2^*([D])$.
Analogously, in degree 12 we have again
$$\bar H_{12}((D\times M)\cup (M\times D);\Q)\cong \bar H_{12}(D\times M;\Q)\oplus\bar H_{12}(M\times D;\Q) \cong \Q\la a_1^*([R]),a_2^*([R])\ra,$$
where the first isomorphism is a consequence of the fact that the two space have the same dimension.

We can construct also a desingularization of $\Sigma$,
$$\tilde\Sigma:= \{(p,v)\in \Pp1\times\Pp1\times V: \text{the curve defined by }v=0 \text{ is singular at } p\}$$
$$\begin{CD}\tilde\Sigma@>\pi>>\Sigma\\
@V\nu VV\\
\;\Pp1\times\Pp1.\end{CD}$$

The map $\pi$ gives $\tilde\Sigma$ the structure of a $\C^{13}$-bundle over $\Pp1\times\Pp1$. Hence $\tilde\Sigma$ is a desingularization of $\Sigma$, and is homotopy equivalent to $\Pp1\times\Pp1$. This ensures that the  
Borel-Moore homology group of degree 28 of $\tilde\Sigma$ is generated by the fundamental classes $[{\tilde D}_1],[{\tilde D}_2]$, where ${\tilde D}_1=\nu^{-1}(\{p_0\}\times\Pp1)$, ${\tilde D}_2=\nu^{-1}(\Pp1\times\{q_0\})$. 

If we choose a point $q_1\in\Pp1\setminus\{q_0\}$, the orbit inclusion defines a map
$$\begin{matrix}\rho_1:&\tilde D&\longrightarrow&\tilde\Sigma\\
&(p,A)&\longmapsto&((p,q_1),\rho (A,I)).\end{matrix}$$

The map $\rho_1$ is well defined, because $\rho(A,I)$ is the union of the line $\{p\}\times\Pp1$, with multiplicity 3, and three lines of the other ruling. Hence it is always singular at $(p,q_1)$.

Since we have 
$$\rho_1^{-1}({\tilde D}_1)= \tau^{-1}(\{p_0\}),$$
$$\rho_1^{-1}({\tilde D}_2)= \emptyset,$$
we know $\rho_1^*([{\tilde D}_1])=[\tau^{-1}(\{p_0\}]$, $\rho_1^*([{\tilde D}_2])=0$. 

Analogously, if we fix $p_1\neq p_0$, the map
$$\begin{matrix}\rho_2:&\tilde D&\longrightarrow&\tilde\Sigma\\
&(p,A)&\longmapsto&((p_1,q),\rho (I, A))
\end{matrix}$$
satisfies $\rho_2^*([{\tilde D}_2])=[\tau^{-1}(\{q_0\}]$, $\rho_2^*([{\tilde D}_1])=0$.

We define $D_1$ and $D_2$ as the images of $\tilde D_1$ and $\tilde D_2$ respectively under the map $\tilde\Sigma\rightarrow\Sigma$.
The following diagrams are commutative:
$$\xymatrix{
{D}\ar[r]^{i_1}&{D\times M}\ar[r]^{\rho}&\Sigma\\
{\tilde D}\ar[u]\ar[rr]^{\rho_1}&&{\tilde\Sigma,}\ar[u]}$$

$$\xymatrix{
{D}\ar[r]^{i_2}&{M\times D}\ar[r]^{\rho}&\Sigma\\
{\tilde D}\ar[u]\ar[rr]^{\rho_2}&&{\tilde\Sigma.}\ar[u]}$$

This implies that we can use our knowledge about $\rho_1$ and $\rho_2$ to study the map induced by $\rho$ on Borel-Moore homology. Thus we find

$$\begin{array}{l@{=}l}
\rho^*([\Sigma])&a_1^*([D])+a_2^*([D]),\\
\rho^*([D_1])&a_1^*([R]),
\\\rho^*([D_2])&a_2^*([R]),
\end{array}$$
which is exactly what we wanted to prove.

\subsection{Application of  Gorinov-Vassiliev's method}\label{quadseq}

We can associate to each $v\in V$ its zero locus on the Segre embedding of ${\Pp1}\times{\Pp1}$, which is a non-singular quadric $Q\subset\Pp3$. 
For $v\neq 0$ we always get a curve of type $(3,3)$ on $Q$. Recall that $Q$ has two rulings, and that curves on $Q$ are classified by considering the number of points of intersection with a general line of respectively the first and the second ruling. 
We want to know what the singular locus of such a curve can be. 
The singular points of a curve are the union of the singular points of each irreducible component of it and the points that are pairwise intersection of components. Note that for each component, the number of singular points allowed is bounded by the aritmetic genus of the component. By writing all possible ways a (3,3)-curve can decompose into components, we get all possibilities on Table \ref{sing33}.

Note that for all configurations of singularities in the same class, the elements of $V$ which are singular at least at the chosen configurations is always a vector space, and has always the same codimension $c$. 
We write this codimension in the second column of the table.

\begin{table}
\caption{\label{sing33}Singular sets of (3,3)-curves}
\centering
\begin{tabular}{|c|c|p{10,5cm}|}
\hline
type& $c$ &description\\\hline
1&3&One point\\\hline
2&6&Two points\\\hline
3&7&Three collinear points\\\hline
4&8&A line\\\hline
5&9&Three non-collinear points\\\hline
6&10&Three collinear points plus a point not lying on the same line\\\hline
7&11&Four points lying on a non-singular conic\\\hline
8&11&{ {5 points: two pairs of points lying each on a line of a different ruling, and the intersection of the two lines}}\\\hline
9&11&A line plus another point\\\hline
10&12&4 points, not lying on the same plane in $\Pp3$\\\hline
11&12&A conic (possibly degenerating into two lines)\\\hline
12&13&5 points: 3 points on a line of some ruling and 2 more points not coplanar with them\\\hline
13&14& {5 points: a point $p$ and 4 more points lying on a non-singular conic not passing through $p$. These points have to be distinct from the intersection of the conic with the two lines of the rulings of $Q$ that pass through $p$}\\\hline
14&14& {6 points: two points on a line of the first ruling, two on a line of the second ruling, the intersection of the two lines and an additional point, not lying on the same line of a ruling of $Q$ as any of the others}\\\hline
15&14& {6 points: union of two triplets of collinear points (on lines of the same ruling). No two points of the configuration can lie on the same line of the other ruling}\\\hline
16&14&A line and two non-collinear points\\\hline
17&15&5 points in general position (no 2 of them on a line of any ruling, no 4 on a conic)\\\hline
18&15& {6 points: two points lying on a line, three points lying on a nonsingular conic and the intersection of the line and the conic. In this configuration only the first two points lie on a line contained in $Q$}\\\hline
19&15&6 points: intersection of $Q$ with three concurrent lines in $\Pp3$\\\hline
20&15& {7 points, intersection of components of the union of two lines of a ruling, one of the other and an irreducible curve of type $(2,1)$ or $(1,2)$}\\\hline
21&15& {7 points, intersection of components of two lines of different rulings and two non-singular conics on $Q$}\\\hline 
22&15& {8 points, intersection of components of the union of two pairs of lines of different rulings and an irreducible conic}\\\hline
23&15&9 points lying on 6 lines (3 in every ruling)\\\hline
24&15&A line and 3 points on a line of the same ruling\\\hline
25&15&A conic (possibly reduced) plus an extra point\\\hline
26&16&   The whole $Q$\\\hline
\end{tabular}
\end{table}

Each item in Table \ref{sing33} can be used to define a family of configurations on $Q$. For every $j=1,\dots,26$, we can define $X_j$ as the space of all configurations of type $(j)$. In this way we get a sequence of subsets $X_1,X_2,\dots,X_{26}$, which satisties conditions 1-4 and 6 in List \ref{ax}. Conditions 5 and 7 do not hold. This problem can be solved by enlarging the list. 

In order for condition \ref{subconf} to be satisfied, it suffices to include all subconfigurations of finite configurations of Table \ref{sing33}.
Condition \ref{border} is more delicate. We have to consider all possible limit positions of configurations of singular points. For instance, points in general position can become collinear, and conics of maximal rank can degenerate to the union of two lines. 

In this way, we can construct a list of configurations that verifies all conditions in List \ref{ax}. This list is really long, so we do not report it here.
Luckily, most configurations give no contribution to the Borel-Moore homology of $\Sigma$. In particular, by Lemma \ref{lem1}, configurations with more than 2 points on a rational curve give no contribution, and, by Lemma \ref{lem2}, the same holds for configurations with at least 2 points on a rational curve minus a point (which is $\cong\C$). Also all configurations containing curves give no contribution, because we can apply Lemma \ref{lemline} and Lemma \ref{2lines}, or a slight modification of them. 

We give below a list limited to the interesting configurations:

\begin{List}
\renewcommand{\labelenumi}{(\Alph{enumi})}
\begin{enumerate}
\item One point.
\item Two points.
\item Three points in general position. 
\item Four points in general position.
\item Six points, intersection of $Q$ with three concurrent lines.
\item Eight points, intersection of components of two pairs of lines for every ruling and a conic. 
\item The whole quadric.

\end{enumerate}
\end{List}
A configuration is said to be \emph{general} if no three points of it lie on the same line, and no four points of it lie on the same plane in $\Pp3$.

\subsection{Columns (A)-(D)}
The space
$F_A$ is a $\C^3$-bundle over $X_A\cong {\Pp1}\times{\Pp1}$. \\
The space $F_B$ is a $\C^6\times\op\Delta_1$-bundle over $X_B\cong \B 2{{\Pp1}\times{\Pp1}}$.\\
The space $F_C$ is a $\C^9\times\op\Delta_2$-bundle over $X_C$, which has the same twisted Borel-Moore homology as $\B 3{{\Pp1}\times{\Pp1}}$ (the non-general configurations form a space with trivial twisted Borel-Moore homology).\\
Analogously, $F_D$ is a $\C^{12}\times\op\Delta_3$-bundle over $X_D$, which has the same twisted Borel-Moore homology as $\B 4{{\Pp1}\times{\Pp1}}$.

Recall that the simplicial bundles are non-oriented, so that we have to consider the Borel-Moore homology with coefficients in the local system $\pm\Q$. This means that we can compute all the terms in these columns by Lemma~\ref{lemQ}.

\subsection{Column (E)}

Configurations of type (E) are the singular loci of (3,3)-curves which are the union of three conics lying on $Q$. 
This gives three pairs of points on the Segre quadric, which are the intersection of it with three concurrent lines in $\Pp3$. Then it is natural to consider this configuration space as a fiber bundle $X_E$ over $\Pp3\setminus Q$. The projection is given by mapping each configuration to the common point of the three lines. 
The fiber over a point $p\in\Pp3\setminus Q$ is the space of configurations of three lines through $p$, all three not tangent to $Q$. It is isomorphic to the space $\tB {3}{\Pp2\setminus C}$ of configurations of three non-collinear points on $\Pp2\cong\Pp{}(T_p\Pp3)$ minus an irreducible conic $C\subset\Pp2$. Note that the orientation of the simplicial bundle over $X_E$ changes under the action of a non-trivial element of the fundamental group of $\Pp2\setminus C$, which is $S_2$. 
 This means that we have to compute the Borel-Moore homology of $\tB 3{\Pp2\setminus C}$ with a 
 system of coefficients locally isomorphic to $\Q$, which changes its orientation when a point of the configuration moves over a loop in $\Pp2\setminus C$ with non-trivial homotopy class. We denote this local system by $W$.

We consider $\tB 3{\Pp2\setminus C}$ as a quotient of the corresponding space of ordered configurations, $\tilde F=(\Pp2\setminus C)^3\setminus \Delta$, where
$$\Delta=\{(x,y,z)\in(\Pp2\setminus C)^3:\dim\la x,y,z\ra\leq1\}.$$

By abuse of notation, we will denote by $W$ also the local system on $(\Pp2\setminus C)^3$.

The closed immersion $\Delta\hookrightarrow(\Pp2\setminus C)^3$ induces the exact sequence in Borel-Moore homology
\begin{equation}\label{hih}
\cdots\rightarrow\bar H_{i+1}(\tilde F;W)\rightarrow \bar H_i(\Delta;W)\rightarrow \bar H_i((\Pp2\setminus C)^3;W)\rightarrow \bar H_i(\tilde F;W)\rightarrow\cdots
\end{equation}
The map associating to a point $p\in(\Pp2\setminus C)$ the intersection with $C$ of the line polar to $p$ with respect to $C$, induces an isomorphism $(\Pp2\setminus C)\cong \B 2C\cong \B 2{\Pp1}$. In particular, we have
\[\bar H_\pu((\Pp2\setminus C)^3;W)\cong (\bar H_\pu(\B 2{\Pp1};\pm \Q))^{\otimes 3},\]
which implies that $\bar H_\pu((\Pp2\setminus C)^3;W)$ is $\Q(3)$ in degree 6, and trivial in all other degrees. Note that the generator of $\bar H_6((\Pp2\setminus C)^3;W)$ is invariant under the action of $S_3$.

We compute next the Borel-Moore homology of $\Delta$ by considering the following filtration:
\begin{equation*}\Delta=\Delta_3\supset\Delta_2\supset\Delta_1,\end{equation*}
$$\Delta_1:=\{(x,y,z)\in\Delta:x=y=z\}$$
$$\Delta_2:=\{(x,y,z)\in\Delta:\exists l\in C\duale(x,y,z\in l)\}.$$

The first term 
$\Delta_1$ of the filtration is isomorphic to $\Pp2\setminus C$, and hence to $\B 2{\Pp1}$, so that its Borel-Moore homology with coefficients in $W$ is $\Q(1)$ in degree 2.

Let us come back to the original construction, and consider ordered configurations of three lines with a common point $p$. Then the configurations in $\Delta_2\setminus \Delta_1$ correspond to triples of lines lying on a plane tangent to $Q$. This means that the space is fibered over the family of such planes, which is parametrized by $Q\duale\cap p\duale$, a non-singular conic. In particular, $Q\duale\cap p\duale$ is simply connected, hence the only system of coefficients we can have there is the constant one.

The fibre over a plane $\Pi$ is given by all ordered triples of lines passing through $p$, lying in $\Pi$ and not tangent to $Q$, such that not all lines coincide. This is the space $(\C^3\setminus \{(z_1,z_2,z_3)\in\C^3:z_1=z_2=z_3\})$, which has Borel-Moore homology $\Q(1)$ in degree 3 and $\Q(3)$ in degree 6.
All elements of $\bar H_\pu(\{z_1=z_2=z_3\};\Q)$ are invariant with respect to the $S_3$-action.

We can conclude that the Borel-Moore homology of $\Delta_2\setminus \Delta_1$ is
$$\bar H_k(\Delta_2\setminus \Delta_1;\Q)=\left\{\begin{array}{ll} \Q(1) & k=3,\\
\Q(2)& k=5,\\
\Q(3)&k=6,\\
\Q(4)&k=8,\\
\,0&\text{else.}\end{array}\right.
$$
For later use, we consider the action induced by the involution $\upsilon$ interchanging the two rulings of $\Pp1\times\Pp1$ on $\bar H_8(\Delta_2\setminus\Delta_1;\Q)$. This Borel-Moore homology group is obtained as the tensor product of $\bar H_4(\Q\duale\cap p\duale;\Q)$ and the Borel-Moore homology group of degree 6 of the fiber. Both factors are invariant under the action induced by $\upsilon$, hence the whole group is invariant under it. Note that the action of $\upsilon$ on the fiber can be seen as interchanging the two lines of intersection of the plane $\Pi$ and $Q$. This means that, if we consider the configuration of points of intersections of the three lines and $Q$, the action of $\upsilon$ interchanges three pairs of points.
 
The space $\Delta_3\setminus \Delta_2$ is a fiber bundle over $\check {\Pp{}}^2\setminus C\duale$. 
The fiber is isomorphic to $(\C^*)^3\setminus \delta$, $\delta:=\{(x,y,z)\in(\C^*)^3:x=y=z\}$. The local system $W'$ induced by $W$ 
coincides with that induced by the local system $T$ on $\C^*$ locally isomorphic to $\Q$ and changing its sign if the point moves along a loop in $\C^*$ whose homotopy class is an odd multiple of the generator of $\pi_1(\C^*)\cong\Z$.

We have again an exact sequence
$$\cdots \!\rightarrow\!\bar H_{i+1}((\C^*)^3\setminus \delta;W')\!\rightarrow\!\bar H_i(\delta,T)\!\rightarrow\!\bar H_i((\C^*)^3,T^{\otimes3})\!\rightarrow\!\bar H_i((\C^*)^3\setminus \delta;W')\,\!\rightarrow\!\,\cdots$$

Since $\bar H_\pu(\C^*,T)$ is trivial by Lemma \ref{lem:T}, $\Delta_3\setminus \Delta_2$ gives no contribution to the cohomology of $\Delta$.

The spectral sequence associated to the filtration $\Delta_i$ is:

\begin{center}\begin{tabular}{c|cc}
6&&$\Q(4)$\\
5&&\\
4&&$\Q(3)$\\
3&&$\Q(2)$\\
2&&\\
1&$\Q(1)$&$\Q(1)$\\
0&&\\
\hline
&1&2
\end{tabular}\end{center}

The only possible non-trivial differential, that from $E_{1,2}$ to $E_{1,1}$, must be an isomorphism for dimensional reasons ($V\setminus \Sigma$ is affine of complex dimension 16). As a consequence,

$$\bar H_k(\Delta;W)^{S_3}=\left\{\begin{array}{ll} \Q(2) &k=5,\\\Q(3)&k=6,\\\Q(4)&k=8,\\0&\text{else.}\end{array}\right.$$

From (\ref{hih}) we get
$$\bar H_k(\tB 3{\Pp2\setminus C};W)=\left\{\begin{array}{ll}\Q(2)&k=6,\\\Q(4)&k=9,\\\,0&\text{else.}\end{array}\right.$$

The Borel-Moore homology of $X_E$ is then by Leray's theorem the tensor product of the above Borel-Moore homology and that of $\Pp3\setminus Q$, which is $\Q(1)$ in degree 3 and $\Q(3)$ in degree 6.
The local system induced on $\Pp3\setminus Q$ is indeed the constant one.

We consider the action induced by the involution $\upsilon$ on the Borel-Moore homology of $X_E$ and $F_E$. We are interested only in determining it on the highest degree component of them.
It is easy to find that $\bar H_3(\Pp3\setminus Q;\Q)$ is anti-invariant with respect to the involution $\upsilon$, and $\bar H_6(\Pp3\setminus Q;\Q)$ is invariant. The group $\bar H_{15}(X_E;\pm\Q)$ is the tensor product of $\bar H_6(\Pp3\setminus Q;\Q)$ and $\bar H_9(\tB 3{\Pp2\setminus C};W)$. The latter group is isomorphic to $\bar H_8(\Delta_2\setminus\Delta_1;\Q)$, which we already showed to be invariant. We also saw that the action of $\upsilon$ on $\Delta\setminus\Delta_1$ produces the interchange of an \emph{odd} number of pairs of points, so that the action of $\upsilon$ on the generator of the Borel-Moore homology of the fiber of the bundle $\Phi_E\rightarrow X_E$ changes its orientation. Hence the invariance of $\bar H_{15}(X_E;\pm\Q)$ under $\upsilon$ implies the anti-invariance of both $\bar H_{20}(\Phi_E;\Q)$ and $\bar H_{22}(F_E;\Q)$.

\subsection{Column (F)}

The configurations of type (F) are the configurations of 8 points which are intersection of components of two pairs of lines for every ruling (say, $l_1,l_2,m_1,m_2$) and a conic $C$. The situation has to be sufficiently general to give exactly 8 points of intersection.

First let us fix the conic. If it is singular, each pair of lines must be an element of $\B 2\C$, which has no Borel-Moore homology with twisted coefficients. 
This means that we have to consider only non-singular conics. On each such conic, the configuration is univocally determined by the intersection points of the 4 lines with it. 
Denote by $\Psi$ the space of two pairs of points in $\Pp1$, i.e., it is the quotient of $\F 4{\Pp1}$ by the relation $(z_1,z_2,z_3,z_4)\sim(z_2,z_1,z_3,z_4)$, $(z_1,z_2,z_3,z_4)\sim(z_1,z_2,z_4,z_3)$. 
The configuration space we have to consider is isomorphic to the product $(\check {\Pp{}}{}^3\setminus Q\duale)\times \Psi$. 

The local system of coefficients is $\Q$ on $\check{\Pp{}}{}^3\setminus Q\duale$ and the rank one local system $R$ on $\Psi$ that changes its sign along the loops exchanging respectively only the first pair of points, or only the second. 

\begin{lem}
The Poincar\'e polynomial of the Borel-Moore homology of $\F 4{\Pp1}$ is $t^4(2+t)(1+t^3)$.
\end{lem}

\proof 
The action of $\PGL(2)$ induces a quotient map $\F 4{\Pp1}\rightarrow \M_{0,4}\cong \Pp1\setminus \{0,1,\infty\}$. 
Notice that there is a natural action of $S_4$ on both $\F 4{\Pp1}$ and $\M_{0,4}$, and that the quotient map is equivariant with respect to this action.
The $S_4$-quotient of $\F4{\Pp1}$ is the space $\Pp{}U_{4,1}$ of non-zero square-free homogeneous quartic polynomials in two variables, up to scalar multiples. 
By \cite{PS} the map $\coh\pu {{\Pp{}} U_{4,1}} \rightarrow \coh\pu{\PGL(2)}$ 
is an isomorphism (it is induced by any orbit map). 
Hence also the pull-back of the orbit map $\coh\pu{\F 4{\Pp1}}\rightarrow\coh\pu{\PGL(2)}$ is surjective.

We can apply the generalized Leray-Hirsch Theorem \ref{ps}, and get that the rational cohomology of $\F 4{\Pp1}$ is isomorphic to $\coh\pu{\M_{0,4}}\otimes_\Q \coh\pu{\PGL(2)}$.
The cohomology of $\PGL (2)$ is $\Q$ in degree 0 and $\Q(-2)$ in degree 3.
Then the claim follows from the fact that $\F 4{\Pp1}$ is smooth of complex dimension 4, so that the cap product with the fundamental class $[\F4{\Pp1}]$ induces an isomorphism $\bar H_\pu(\F 4{\Pp1};\Q)\cong\coh{8-\pu}{\F 4{\Pp1}}(4)$.
\qed

\begin{lem}\label{qqq}
$\bar H_\pu(\Psi;R)$ is $\Q(1)$ in degree 4 and $\Q(3)$ in degree 7. 
Moreover, $\bar H_\pu(\Psi;R)$ is invariant with respect to the $S_2$-action induced by the interchange of $(z_1,z_2,z_3,z_4)$ and $(z_3,z_4,z_1,z_2)$.
\end{lem}

\proof

We refer to the preceding Lemma. 

$\bar H_\pu(\Psi;R)$ can be identified with the part of $\bar H_\pu(\F 4{\Pp1};\Q)$ which is anti-in\-var\-iant with respect to the action of the transpositions $(1,2)$ and $(3,4)$. By the equivariance of the quotient map, it is sufficient to determine the action on $\M_{0,4}$.
The action of $S_4$ on $\coh\pu{\M_{0,4}}$ factorizes via $S_3$, hence the action of any pair of commuting transpositions must coincide.
The anti-invariant part of $\coh\pu{\M_{0,4}}$ with respect to $(1,2)\in S_4$ is $\Q(-1)$ in degree 1. 
The actions of the two transpositions $(1,3)$ and $(2,4)$ on $\M_{0,4}$ coincide, and their product is the identity.
If we pass from cohomology to Borel-Moore homology, this implies the claim.
\qed

Lemma \ref{qqq} implies that the Borel-Moore homology of $X_{F}$ is
$$\bar H_k(X_F;\pm\Q)=\left\{\begin{array}{ll}\Q(2)&k=7,\\\Q(4)^2&k=10,\\\Q(6)&k=13,\\\,0&\text{else.}\end{array}\right.$$

Moreover, we can compute the action on $\bar H_\pu(X_F;\pm\Q)$ induced by the involution $\upsilon$ interchanging the two rulings of the quadric $Q$.
Fix a configuration in $X_F$. Up to a choice of coordinates, we can assume that the action of $\upsilon$ permutes the two lines of the first ruling with the two of the second ruling in the configuration. This means that the action of $\upsilon$ on $\Psi$ can be identified with that interchanging $(z_1,z_2,z_3,z_4)\in \F4{\Pp1}$ with $(z_3,z_4,z_1,z_2)$. As a consequence, $\bar H_{13}(X_F;\pm\Q)$ is invariant with respect to $\upsilon$. If we look at the configuration itself, three pairs of points are exchanged. Namely, the four points of intersection of the conic and the lines are exchanged in pairs, and also exactly two of the points of intersection between the lines are exchanged. Hence the Borel-Moore homology of the fiber of the bundle $\Phi_F\rightarrow X_F$ is anti-invariant for $\upsilon$. We can conclude that the highest class in the Borel-Moore homology of $F_F$ (i.e., that of degree 22) is anti-invariant for the action of $\upsilon$.

\subsection{Column (G)}

By Proposition \ref{ucci}, the space $F_G$ is an open cone. The Borel-Moore homology of its base space can be computed by the spectral sequence in Table~\ref{pic1}. Its columns coincide with those of the main spectral sequence, but are shifted by twice the dimension of the fibre of the complex vector bundle we considered for each of them.

\begin{table}
\caption{\label{pic1}Spectral sequence converging to the base of $F_G$.}
\centering
\begin{tabular}{r|p{24pt}|p{24pt}|p{24pt}|p{24pt}|p{24pt}|p{24pt}}
15&&&&&$\Q(7)$\\\hline
14&&&&&&$\Q(6)$\\\hline
13&&&&&\\\hline
12&&&&&$\Q(5)^2$&\\\hline
11&&&&&&$\Q(4)^2$\\\hline
10&&&&&&\\\hline
9&&&&&$\Q(3)$&\\\hline
8&&&&&&$\Q(2)$\\\hline
7&&&$\Q(4)$&$\Q(4)$&\\\hline
6&&&&&\\\hline
5&&$\Q(3)^2$&$\Q(3)^2$&&\\\hline
4&&&&&\\\hline
3&$\Q(2)$&$\Q(2)^2$&$\Q(2)$&&\\\hline
2&&&&&\\\hline
1&$\Q(1)^2$&$\Q(1)^2$&&&\\\hline
0&&&&&\\\hline
-1&$\ \,\,\Q$&&&&\\
\hline
&\ \,A&\ \,\,B&\ \,\,C&\ \,\,D&\ \,\,E&\ \,\,F
\end{tabular}
\end{table}

For dimensional reasons, the rows of indices 1, 3, 5 and 7 of this exact sequence are exact. The reason is that if were not so, there would be non-trivial elements in the general spectral sequence, in a position such that they could not disappear. Hence they would give non-trivial elements in $H^{j}(V\setminus \Sigma)\cong\bar H_{31-j} (\Sigma)$ for $j>16$, which is impossible because $V\setminus \Sigma$ is affine of dimension 16.

As a consequence, the spectral sequence here degenerates at $E^2$. 
By Lemma~\ref{opencone}, the Borel-Moore homology of the open cone is

$$\bar H_k(F_G;\Q)=\left\{\begin{array}{ll}\Q(2)+\Q(3)&k=15,\\\Q(4)^2+\Q(5)^2&k=18,\\\Q(6)+\Q(7)&k=21,\\0&\text{else.}\end{array}\right.$$

\subsection{Spectral sequence}

\begin{table}
\caption{\label{gr1} Spectral sequence converging to $\bar H_k(\Sigma;\Q)$.}
\centering
\begin{tabular}{r|p{28pt}|p{28pt}|p{28pt}|p{28pt}|p{28pt}|p{28pt}|c}
29&$\Q(15)$&&&&&&\\\hline
28&&&&&&&\\\hline
27&$\Q(14)^2$&&&&&&\\\hline
26&&&&&&&\\\hline
25&$\Q(13)$&$\Q(13)^2$&&&&&\\\hline
24&&&&&&&\\\hline
23&&$\Q(12)^2$&&&&&\\\hline
22&&&&&&&\\\hline
21&&$\Q(11)^2$&$\Q(11)$&&&&\\\hline
20&&&&&&&\\\hline
19&&&$\Q(10)^2$&&&&\\\hline
18&&&&&&&\\\hline
17&&&$\Q(9)$&&$\Q(8)$&&\\\hline
16&&&&&&$\Q(7)$&\\\hline
15&&&&$\Q(8)$&&&\\\hline
14&&&&&$\Q(6)^2$&&$\Q(7)+\Q(6)$\\\hline
13&&&&&&$\Q(5)^2$&\\\hline
12&&&&&&&\\\hline
11&&&&&$\Q(4)$&&$\Q(5)^2+\Q(4)^2$\\\hline
10&&&&&&$\Q(3)$&\\\hline
9&&&&&&& \\\hline
8&&&&&&&$\Q(3)+\Q(2)$ \\\hline
&\ \,\,\,A&\ \,\,\,\,B&\ \,\,\,\,C&\ \,\,\,\,D&\ \,\,\,\,E&\ \,\,\,\,F&G
\end{tabular}
\end{table}

We know now all columns of the spectral sequence associated to the filtration $$\Fil_A(\ba\X(\pu)\ba)\subset\cdots\subset\Fil_G(\ba\X(\pu)\ba).$$ Its $E^1$ term is represented in Table \ref{gr1}.

This spectral sequence degenerates at $E^1$. Indeed, the only possible non-trivial differentials are between the first 4 columns. We know a priori that $\coh 0 {C_0}$ has dimension 1. Then Theorem \ref{ps} implies that the cohomology of $X=V\setminus \Sigma$ must contain a copy of the cohomology of the group $G$. This is impossible if any of the differentials in columns (A)-(D) is non-zero.

\begin{table}
\caption{\label{coho} Rational cohomology of $X$ and $G$.}
$$\begin{array}{l@{(X;\Q)=}ll@{(G;\Q)=}l}
H^0& \Q&H^0&\Q\\
H^1& \Q(-1)&H^1&\Q(-1)\\
H^2 & \,0& H^2 & \,0\\ 
H^3& \Q(-2)^2& H^3&\Q(-2)^2\\
H^4& \Q(-3)^2& H^4&\Q(-3)^2\\
H^5& \Q(-3)&H^5&\,0\\ 
H^6& \Q(-4)^2&H^6&\Q(-4)\\
H^7& \Q(-5)& H^7&\Q(-5)\\
H^8& \Q(-5)^2 \\
H^9& \Q(-6)^2+\Q(-8)+\Q(-9)\\
H^{10}& \Q(-9)+\Q(-10)\\
H^{11}& \Q(-7)\\
H^{12}& \Q(-8)+\Q(-10)^2+\Q(-12)^2\\
H^{13}& \Q(-11)^2+\Q(-12)^2\\
H^{14}&\,0\\
H^{15}& \Q(-12)+\Q(-11)\\
H^{16}& \Q(-12)+\Q(-13)
\end{array}$$
\end{table}

We can compute the whole cohomology of $X$ from the Borel-Moore homology of $\Sigma$, 
using the isomorphism induced by the cap product with the fundamental class of the discriminant
$$\tilde H^\pu (X;\Q)\cong \bar H_{31-\pu}(\Sigma;\Q)(-d).$$
The results are in Table~\ref{coho}. 
By Theorem \ref{ps}, the rational cohomology of $X/G$ is
$$\begin{array}{l}
\coh 0{X/G}=\Q\\
\coh 5{X/G}=\Q(-3)\\
\coh 9{X/G}=\Q(-8)+\Q(-9)
\end{array}$$
and trivial in all other degrees. 

The quotient $X/G$ is a double cover of $C_0$, the $S_2$-action being generated  by the involution $\upsilon$ interchanging the two rulings of the Segre quadric $Q$. The cohomology of $X$ is invariant with respect to this involution in the degrees $0,5$.
The cohomology in degree 9 comes from the terms $E^1_{E,17}$, $E^1_{F,16}$ in the spectral sequence.
During the computation of columns (E) and (F) we observed that these terms are anti-invariant with respect to the action of $S_2$ on $X$. Then $C_0$ has no cohomology in degree 9.

We can conclude that the rational cohomology of $C_0$ is $\Q$ in degree 0 and $\Q(-3)$ in degree 5, as we claimed in Theorem \ref{c0}.

\section{Curves on a quadric cone}\label{C_1}

\subsection{The space $C_1$ as geometric quotient}
The aim of this section is to realize $C_1$ as a geometric quotient, satisfying the hypotheses of Theorem \ref{ps}. 
We perform it by regarding the elements of $C_1$ as curves of degree 6 on a quadric cone.

After a choice of coordinates, we can identify the quadric cone with $\Pp{}(1,1,2)$. Then we consider the vector space $\C[x,y,z]_6$, where $\deg x=\deg y=1$, $\deg z=2$. 

This space has complex dimension 16. 
The polynomials defining singular curves form the \emph{discriminant} $\Sigma$, which has in this case two irreducible components of dimension 15. One component, which we denotes by $H$, is the locus of curves passing through the vertex of the cone. Such a curve is always singular, and in general this singularity can be solved by considering the proper transform in the blowing up of the vertex of the cone.
The other component, which we denote by $S$, is the closure in $\C^{16}$ of the locus of curves which are singular in a point different from the vertex of the cone. More specifically, $S$ is the locus of curves such that the proper transform in the blowing up of the vertex of the cone is singular, or tangent to the exceptional locus.
Note that $H$ is the hyperplane defined by the condition that the coefficient of $z^3$ is zero. The other component $S$ is the affine cone over the dual variety of the Veronese embedding of $\Pp{} (1,1,2)$ in $\Pp{15}$.

We find next generators of the cohomology of $X$ in degree 0 and 2, which have an interpretation as fundamental classes of subvarieties of $S\cup H$. The whole cohomology of $X=\C^{16}\setminus (S\cup H)$ will be calculated in Section \ref{conseq}. 

We can compute the Borel-Moore homology of $S\cup H$ by considering the incidence correspondence 
$$\mathcal T=\{(f,p)\in (S\cup H)\times\Pp{}(1,1,2):f \text{ is singular at }p\}.$$

The \  incidence \  correspondence \  has \  a \  natural \  projection \  $\pi$ \  to \  the \  cone 
$\Pp{}(1,1,2)$. The fiber over the vertex $[0,0,1]$ is simply $\{[0,0,1]\}\times H$. If we restrict to $\pi^{-1}\left(\Pp{}(1,1,2)\setminus \{[0,0,1]\}\right)$, $\pi$ is a complex vector bundle of rank 13. This allows us to compute the Borel-Moore homology of $\mathcal T$:
$$\bar H_k(\mathcal T;\Q)=\left\{\begin{array}{ll}\Q(15)^2 &k=30\\\Q(14)&k=28\\0&\text{else.}\end{array}\right.$$

Although $\pi$ is no desingularization of $S\cup H$, it is a proper map, hence it induces a homomorphism of Borel-Moore homology groups. By the results in Section \ref{conseq} 
the map induced by $\pi$ in those degrees must be an isomorphism.
In particular, this implies that $\bar H_{28}(S\cup H;\Q)$ is generated by the fundamental class of the locus $S'$ of curves with a singularity on a chosen line. Obviously $\bar H_{30}(S\cup H;\Q)$ is generated by $[S]$ and $[H]$.

Automorphisms of the graded ring $\C[x,y,z]$ are of the form
$$\left\{\mil\begin{array}{cll}x&\mapsto&\alpha x +\beta y\\
y & \mapsto &\gamma x + \delta y\\
z & \mapsto &\epsilon z + q(x,y),
\end{array}\right.$$
where $\alpha,\beta,\gamma,\delta,\epsilon$ are complex numbers such that $\epsilon(\alpha\delta-\beta\gamma)\neq0$ and $q\in\C[x,y]_2$. 

They form a group $G$ of dimension 8.
By contracting the vector space $\C[x,y]_2\cong\C^3$ to a point, we get that $G$ is homotopy equivalent to $\GL(2)\times\C^*$. This means that we can apply Theorem \ref{ps} to the action of $\GL(2)\times\C^*$ instead of the whole $G$. 
In order to be able to apply the Theorem in Section~\ref{conseq}, we check that its hypothesis are satisfied. Namely, for each generator $\eta$ of degree $2r-1$ of the cohomology of $\GL(2)\times\C^*$, we want to define a subscheme $Y$ of the discriminant, of pure codimension $r$, whose fundamental class maps to a non-zero multiple of $\eta$ under the composition
$$\bar H_{2(16-r)}(Y)\rightarrow \bar H_{2(16-r)}(\Sigma) \xrightarrow{\sim} H^{2r-1}(X)\xrightarrow{\rho^*}H^{2r-1}(G),$$
where $\rho$ denotes any orbit inclusion of $G$ in $\C^{16}$.

The cohomology of the product $\GL(2)\times\C^*$ decomposes naturally into that of its subgroups $\{I\}\times\C^*$ and $\GL(2)\times\{1\}$. We can reduce the study of the orbit inclusion to that of the two maps
$$\begin{matrix} \rho_1: &\C^*&\longrightarrow&X\\
& t&\longmapsto&\rho(I,t)\end{matrix}$$
and
$$\begin{matrix} \rho_2: &\GL(2)&\longrightarrow&X\\
& A&\longmapsto&\rho(A,1).\end{matrix}$$

The map induced by $\rho_1$ on cohomology is
$${\bar H_{31-\pu}(H\cup S)}\cong{\coh\pu X}\xrightarrow{\rho_1^*}{\coh\pu {\C^*}.} $$

The cohomology of $\C^*$ is generated by the fundamental class of $\{0\}\subset\C$. If we extend $\rho_1$ to a map $\C\longrightarrow \C^{16}$, we find that 0 is mapped to a curve which is the union of 6 distinct lines through the vertex of the cone. This means that the preimage of $H$ coincides with 0, while the preimage  of $S$ is empty.
Hence $\rho_1^*([H])$ is a non-zero multiple of $[0]$ (in fact, by direct computation, it is $3[0]$).

We consider next $\rho_2$. It induces the map
$${\bar H_{31-\pu}(H\cup S)}\cong{\coh\pu X}\xrightarrow{\rho_2^*}{\coh\pu {\GL(2)}}.$$

Recall from Section \ref{LH_0} that the cohomology of $\GL(2)$ is generated by the fundamental class of the complement $D$ of $\GL(2)$ in the space $M$ of $2\times2$ matrices and that of the subspace $R$ of matrices with only zeroes on the first column. The Borel-Moore homology class $[D]\in \bar H_6(D;\Q)$ corresponds to a class of degree 1 in the cohomology of $\GL(2)$, and the Borel-Moore homology class $[R]$ corresponds to a class in $\coh 3{\GL(2)}$.

We look at the extension of $\rho_2$ to $M\longrightarrow \C^{16}$. 
The elements in $D$ are mapped to curves which are the union of three non-singular quadrics, having the same tangent line in a common point. These are always elements of $S\setminus H$.
If we choose $S'$ as the locus of curves singular at some point of the line $\{y=0\}\subset\Pp{}(1,1,2)$, then we have that the preimage of $S'$ is exactly $R$. These considerations imply the surjectivity of the orbit inclusion on cohomology, hence the hypotheses of Theorem \ref{ps} are established.

\subsection{Application of Vassiliev-Gorinov's method}\label{conseq}

As a starting point for the application of Gorinov-Vassiliev's method, we list the possible singular loci of elements of $\Sigma$. This is achieved by considering all possible decompositions in irreducible components of a curve of degree 6 on $\Pp{}(1,1,2)$. Then we consider how many singular points can lie on each components, and where the pairwise intersection of components can lie. The results are listed in Table~\ref{singcon}.

For two configurations of the same type, the linear subspace of $V$ of curves that are singular in the one and the linear subspace of curves that are singular in the other have always the same codimension $c$. Hence this codimension depends only on the type of the configuration. We write it in the second column of Table~\ref{singcon}.

\begin{table}
\caption{\label{singcon}Singular configurations of degree 6 curves on a quadric cone}
\centering
\begin{tabular}{|c|c|p{10,5cm}|}
\hline
type&$c$&description\\\hline
1&1&The vertex\\\hline
2&3&A general point\\\hline
3&4&The vertex and a general point\\\hline
4&6&The vertex and two points on a line of the ruling\\\hline
5&6&Two general points\\\hline
6&7&A singular line\\\hline
7&7&The vertex and 2 general points\\\hline
8&9&Three general points\\\hline
9&9&The vertex, two points on a line of the ruling and a general point\\\hline
10&10&The vertex and 3 general points\\\hline
11&10&A singular line and a general point\\\hline
12&11&The vertex, two points on a line of the ruling, and two points on another\\\hline
13&11&4 points on a non-singular conic\\\hline
14&12&The union of a line and two points on another line of the ruling\\\hline
15&12&Four general points\\\hline
16&12&The vertex, two collinear points and two general points\\\hline
17&12&The vertex and 4 points on a non-singular conic\\\hline
18&12&Two lines\\\hline
19&12&A non-singular conic\\\hline
20&13&The vertex and a non-singular conic\\\hline
21&13&The union of a line and two general points\\\hline
22&14&The vertex, 4 points on a conic and a point lying on the line joining one of the 4 points with the vertex\\\hline
23&14&The vertex, two points on a line of the ruling, two points on another and a general point\\\hline
24&14&4 points on a non-singular conic and a general point\\\hline
25&15&A rational normal cubic\\\hline
26&15&The intersection of components of the union of two lines and two conics\\\hline
27&15&6 points, intersection of the cone with three concurrent lines\\\hline
28&15&Three lines\\\hline
29&15&Union of a line and a non-singular conic\\\hline
30&16&The whole cone\\\hline
\end{tabular}
\end{table}

In Table~\ref{singcon}, a configuration of points is said to be \emph{general} if it does not contain the vertex, no two points of the configuration lie on the same  line of the ruling of the cone and at most 3 points lie on the same conic contained in the cone. 

The families on configurations defined by Table~\ref{singcon} does not satisfy the conditions in List \ref{ax}. We have to refine it in order to have all finite subconfigurations of finite configurations included.
In this way we get a new list as follows:

\begin{List}
Sequence of families of configurations which satisfies the conditions on List~\ref{ax}. The number under square brackets is the codimension of the space of polynomials singular at a chosen configuration of that type.

\begin{itemize}
\item[(A)]
The vertex.
\ [1]

\item[(B)]
A general point.
\ [3]

\item[(C)]
The vertex and a general point.
\ [4]

\item
Two points distinct from the vertex, lying on the same line of the ruling.
\ [6]

\item
The vertex and two points in the same line of the ruling.
\ [6]

\item[(D)]
Two general points.
\ [6]

\item
Three collinear points, distinct from the vertex.
\ [7]

\item
The vertex and three collinear points.
\ [7]

\item
Four collinear points, distinct from the vertex .
\ [7]

\item
The vertex and four collinear points.
\ [7]

\item
Five collinear points, distinct from the vertex .
\ [7]

\item
The vertex and five collinear points.
\ [7]

\item
A singular line.
\ [7]

\item[(E)]
The vertex and 2 general points.
\ [7]

\item
Two points on the same line of the ruling and a general point.
\ [9]

\item
Three general points.
\ [9]

\item
The vertex, two points on a line of the ruling and a general point.
\ [9]

\item
The vertex and 3 general points.
\ [10]

\item
Three collinear points (distinct from the vertex) and a general point.
\ [10]

\item
The vertex, three collinear points and a general point.
\ [10]

\item
Four collinear points (distinct from the vertex) and a general point.
\ [10]

\item
The vertex, four collinear points and a general point.
\ [10]

\item
A singular line and a general point.
\ [10]

\item
Two points on a line of the ruling, and two points on another.
\ [11]

\item
The vertex, two points on a line of the ruling, and two points on another.
\ [11]

\item
Four points on a non-singular conic.
\ [11]

\item
Two points on a line of the ruling and three points on another. All points are different from the vertex.
\ [12]

\item
The vertex, two points on a line and four points on another.
\ [12]

\item
Two points on a line of the ruling and three points on another. All points are different from the vertex.
\ [12]

\item
The union of a line and 2 points on a line of the ruling.
\ [12]

\item
Two points on a line of the ruling and 2 general points.
\ [12]

\item
Four general points.
\ [12]

\item
The vertex, two collinear points and 2 general points.
\ [12]

\item
The vertex and 4 points on a non-singular conic.
\ [12]

\item
Two lines.
\ [12]

\item
A non-singular conic.
\ [12]

\item
The vertex and a non-singular conic.
\ [13]

\item
The union of a line and two general points.
\ [13]

\item
Four points on a conic and a point lying on the line joining one of the 4 points with the vertex.
\ [14]

\item
The vertex, 4 points on a conic and a point lying on the line joining one of the 4 points with the vertex.
\ [14]

\item
Two points on a line of the ruling, two points on another and a general point.
\ [14]

\item
The vertex, two points on a line of the ruling, two points on another and a general point.
\ [14]

\item
Four points on a non-singular conic and a general point.
\ [14]

\item
A rational normal cubic.
\ [15]

\item[(F)]
The intersection of components of the union of two lines and two conics, excluding the vertex.
\ [15]

\item[(G)]
The intersection of components of the union of two lines and two conics.
\ [15]

\item[(H)]
6 points, intersection of the cone with three concurrent lines.
\ [15]

\item
Three lines.
\ [15]

\item
Union of a line and a non-singular conic.
\ [15]

\item[(I)]
The whole cone.
\ [16]

\end{itemize}
\end{List}

By the results in Sections \ref{lemmas} and \ref{simplices}, or an adaptation of them, the only configurations giving non-trivial Borel-Moore homology are those indicated with (A)-(I). We will consider only them.

\begin{description}
\item[(A)]
Obviously, $F_A\cong\C^{15}$. The only homology is in degree 30.

\item[(B)]
The space $F_B$ is a $\C^{13}$-bundle over the cone minus its vertex. 

\item[(C)]
The space $F_C$ is a $\C^{10}\times{\op\Delta_1}$-bundle over the cone minus its vertex.

\item[(D)]
The space $\Phi_D$ is a non-orientable bundle of open simplices of dimension 1 over the subspace of $\B 2{\C\times{\Pp1}}$ consisting of points not on the same line. This has the same Borel-Moore homology with twisted coefficients of $\B 2{\C\times{\Pp1}}$, which is non-trivial only in degree 6. 
The space $F_D$ is a $\C^{10}$-bundle over $\Phi_D$.

\item[(E)]
The space $\Phi_E$ is a non-orientable bundle of open simplices of dimension 2 over the subspace of $\B 2{\C\times{\Pp1}}$ consisting of points not on the same line. 
The space $F_E$ is a $\C^{9}$-bundle over $\Phi_E$.

\item[(F) and (G)]
We can consider configurations of types (F) and (G) together. We have then that $\Phi_F\cup\Phi_G$ is a complex vector bundle of rank one over $F_F\cup F_G$.

We claim that Borel-Moore homology of $F_F\cup F_G$ is trivial. Let us consider the fiber $\Psi$ of its projection to $X_F$ (which is canonically isomorphic to $X_G$). Then $\Psi$ is a simplex with vertices $t,a_1,a_2,a_3,a_4,a_5,a_6$, where $\{a_i\}\in X_F$ and $t$ is the top of the cone. We have to consider the closed simplex $S$, minus all external faces containing $t$. Denote by $B$ the union of such external faces. Observe that both $B$ and $S$ can be contracted to the vertex $t$. Then $\bar H_\pu(S;\Q)=H_\pu(C,B)=0$, which yields the claim.

\item[(H)]
The space $F_H$ is a $\C\times\op\Delta_5$-bundle over $X_H$. We claim that $X_H$ has no Borel-Moore homology in the system of coefficients induced by $F_H$. 

Each configuration in $X_H$ is determined by three concurrent lines in $\Pp3$. This implies that $X_H$ is fibred over the complement of the cone in $\Pp3$. The fiber over a point $p$ is the space of configurations of three distinct lines through $p$, subject to certain conditions. Namely, the lines cannot be tangent to the cone, any such triple must generate $\Pp 3$, and no two of the lines can generate a plane passing through the vertex of the cone. The condition of being tangent to the cone defines the union of two lines $l,m$ inside $\Pp2\cong\Pp{}(T_p\Pp3)$; then what we need is to compute the Borel-Moore homology of the subset $\tB 3{\C\times\C^*}$ of configurations of triples of points on $\C\times\C^*$, such that the three points are not collinear, and no line passing through two of them contains the point $l\cap m$. 

The system of coefficients we have to consider changes its sign every time that a point moves along a loop around one of the lines $l$ or $m$. 
This system of coefficients is induced by a system $T'$ on
$\C\times\C^*$, and Lemma \ref{lem:T} implies that $\bar H_\pu(\C\times\C^*;T')$ is trivial. Hence also $\B 3{\C\times\C^*}$ has no Borel-Moore homology in the system of coefficient we are interested in.
Thus all we need is to prove that $S=\B 3{\C\times\C^*}\setminus \tB 3{\C\times\C^*}$ has also trivial Borel-Moore homology in the chosen system of coefficients.

We can distinguish three subsets of $S$:\begin{enumerate}
\item Configurations of three points on a line passing through $l\cap m$.\\
 Such lines are parameterized by $\C^*$, hence this subspace of configurations is isomorphic to the product $\C^*\times\B 3\C$. The system of coefficients induced on $\C^*$ is $T$. Recall from Lemma \ref{lem:T} that $\bar H_\pu(\C^*,T)=0$. We can conclude that this kind of configurations give no contribution to the Borel-Moore homology of $S$. 

\item Configurations with three collinear points on a line not passing through $l\cap m$.\\
 This subspace is fibred over $\C^2$ (the space of such lines), the fiber being $\B 3{\C^*}$. The claim follows from the fact that the system of coefficients on $\B 3{\C^*}$ coincides with that induced by $T$.

\item Configurations  with two points on a line passing through $l\cap m$, and a third point outside this line.\\ 
We know that the conditions that the third point is outside the line is not relevant for the Borel-Moore homology of the configuration. Without that condition, this subset of configurations is a product $\B2\C\times\C^*\times (\C\times\C^*)$, where the system of coefficients induced on the factor $\C^*$ is $T$. Then the Borel-Moore of this substratum is trivial.
\end{enumerate}

\item[(I)]
By point (\ref{opeco}) in Proposition \ref{ucci}, the space $F_I$ is an open cone over a space, whose Borel-Moore homology can be computed by the following spectral sequence:

\begin{center}\begin{tabular}{r|c|c|c|c|c|c|c}
 4&    &    &    &    &    &    \\\hline 
 3&    &    &    &$\Q(3)$&$\Q(3)$&    \\\hline 
 2&    &$\Q(2)$&$\Q(2)$&    &    &    \\\hline 
 1&    &    &    &    &    &    \\\hline 
 0&    &$\Q(1)$&$\Q(1)$&    &    &        \\\hline 
-1&$\Q$&    &    &    &    &           \\\hline 
  & A  & B  & C  & D  & E  &  F,G & H  
\end{tabular}\end{center}

Note that all differentials of this spectral sequence must be 0 for dimensional reasons.

\end{description}

Putting information over all columns together, we have that the spectral sequence converging to Borel-Moore homology of $\Sigma=S\cup H$ has the following $E^1$ term:

\begin{center}\begin{tabular}{r|c|c|c|c|c}
30&    &    &    &    &    
\\\hline
29&$\Q(15)$&    &    &    &    
\\\hline
28&    &$\Q(15)$&    &    &    
\\\hline
27&    &    &    &    &    
\\\hline
26&    &$\Q(14)$&$\Q(14)$&    &    
\\\hline
25&    &    &    &    &    
\\\hline
24&    &    &$\Q(13)$&    &    
\\\hline
23&    &    &    &$\Q(13)$&    
\\\hline
22&    &    &    &    &    
\\\hline
21&    &    &    &    &$\Q(12)$
\\\hline
  & A  & B  & C  & D  & E  
\end{tabular}\end{center}

We have omitted columns from (F) to (I) because they give trivial contributions.

The differential $E^1_{C,26}\rightarrow E^1_{B,26}$ is trivial because $\coh \pu X$ must contain a copy of $\coh \pu G$.
This means that this spectral sequence degenerates at $E^1$. 
As a consequence, $\coh \pu X=\coh \pu G$. 
Theorem \ref{ps} gives that $C_1=X/G$ has the rational cohomology of a point. Theorem \ref{c1} is now established.

\section{Hyperelliptic locus}\label{C_2}

The moduli space $\mathcal H_g$ of smooth hyperelliptic curves of genus $g\geq2$ 
 has always the cohomology of a point. We have that $\mathcal H_g$ coincides with the moduli space of configurations of $2g+2$ distinct points on $\Pp1$, which is in turn the quotient of the space of
binary polynomials of degree $2g+2$ without double roots, for the action of $\GL(2)$. 
This is a special case of moduli space of hypersurfaces. It is shown in \cite{PS} that in the hypersurface case the hypotheses of Theorem \ref{ps} are always satisfied. 
Thus all we need to prove is that the cohomology of the space of polynomials without double roots is generated by the elements mapped to the generators of the cohomology of $\GL(2)$. 

\begin{lem}
Let $V=\C[x,y]_d$ be the vector space of homogeneous binary polynomials of degree $d\geq 4$ , and $\Delta\subset V$ the discriminant, i.e., the locus of polynomials with multiple roots. Then the cohomology of $V\setminus \Delta$ with coefficients in $\Q$ is $\Q$ in degree 0, $\Q(-1)$ in degree 1, $\Q(-2)$ in degree 3 and $\Q(-3)$ in degree 4.
\end{lem}

\proof We apply Gorinov-Vassiliev's method and compute the Borel-Moore homology of $\Delta$. For every $v\in \Delta$, we denote by $K_v$ the locus in $\Pp1$ which is the projectivization of the multiple roots of $v$. Denote by $X_k$ the family of all configurations of $k$ distinct points in $\Pp1$. $X_1,\dots,X_{[d/2]},\{\Pp1\}$ satisfy the conditions in List \ref{ax}. The dimension $d_i$ of $L(x)$ for $x\in X_i$ is always $d+1-2i$.
This means that we can construct a geometrical resolution of $\Delta$. Note that in the filtration only the first two terms give non-trivial Borel-Moore homology, because $X_k\cong\B k{\Pp1}$ has trivial twisted Borel-Moore homology for $k\geq3$. 

The spectral sequence for the Borel-Moore homology of $\Delta$ is
$$\begin{array}{l|c|c}
2d-1&\Q(d)&\\\hline
2d-2&&\\\hline
2d-3&\Q(d-1)&\\\hline
2d-4&&\\\hline
2d-5&&\Q(d-2)\\\hline
&1&2
\end{array}$$

Clearly the spectral sequence degenerates at $E^1$. This implies the claim.\qed

\subsection*{Acknowledgements}

I would like to thank my advisor Joseph Steenbrink for proposing the subject and for useful discussions and suggestions.

\end{document}